\documentclass[12pt,reqno,a4paper]{amsart}
\usepackage{extsizes}
\usepackage[a4paper,  top=0.5in, bottom=0.5in, left=0.51in, right=0.51in]{geometry}

\calclayout
\usepackage{graphicx}
\usepackage[utf8]{inputenc}
\usepackage{nicefrac}
\usepackage{mathabx}
\usepackage{enumitem}
\usepackage{subcaption}
\usepackage{amsrefs}
\usepackage{esint}
\usepackage{mathtools}

\usepackage{wrapfig}
\usepackage{xcolor}
\usepackage{cleveref}
\usepackage{enumitem}
	
\usepackage{ulem}
\usepackage{autonum}

\newtheorem{theorem}{Theorem}[section]
\newtheorem{lemma}[theorem]{Lemma}
\newtheorem{corollary}[theorem]{Corollary} 

\theoremstyle{definition}

\newtheorem{prop}[theorem]{Proposition}
\newtheorem{cor}[theorem]{Corollary}

\allowdisplaybreaks
\theoremstyle{remark}
\newtheorem{remark}[theorem]{Remark}

\numberwithin{equation}{section}

\allowdisplaybreaks

\usepackage{amsmath}
\usepackage{braket} 

\makeatletter
\@namedef{subjclassname@2020}{%
  \textup{2020} Mathematics Subject Classification}
\makeatother

\begin{document}

\title[Alt-Caffarelli problems for elastic graphs]{An Alt-Caffarelli problem for elastic graphs}


\author{Marius Müller}
\address[M.~M\"uller]{Universität Augsburg, Institut für Mathematik, Universitätsstraße 14, 86159 Augsburg}
\email{marius1.mueller@uni-a.de}


\date{\today}


\keywords{Alt-Caffarelli Problem, Elastic energy, Free boundary,  Noncoercive Variational Problem. \\ \textit{AI declaration.} During the preparation of this manuscript ChatGPT (OpenAI) was used to check grammar and intermediate computations. The author critically reviewed the suggestions and is fully responsible for mathematical arguments.}

\begin{abstract}
We study an Alt-Caffarelli-type free boundary problem for symmetric elastic graphs with fixed boundary height. Our results concern existence of minimizers and their description in terms of regularity and shape. It turns out that the weight coefficient in front of the measure term is decisive for existence. We identify a sharp threshold for this coefficient below which existence is guaranteed.  Moreover, we find that all sign-changing minimizers are described by a system of three Eulerian elasticae, none of which are lines. Another result shows non-flatness of free boundary points, a classical question in the study of Alt-Caffarelli problems. These properties lead to a better understanding of the nonexistence phenomenon above the aforementioned weight threshold.
\end{abstract}

\maketitle

\section{Introduction}
In this article we seek to minimize an \textit{elastic Alt-Caffarelli} energy 
\begin{equation}
    \mathcal{E}_\theta(u)  := \int_{\mathrm{graph}(u)} \kappa^2 \; \mathrm{d}\mathbf{s} +  \theta |\{ u> 0 \}|,
\end{equation}
among all symmetric graphs $u:(0,1) \to \mathbb{R}$ with a fixed boundary value $h> 0$, i.e. the admissible set is given by
\begin{equation}
    \mathcal{A}_{\rm sym}(h) := \{ u \in W^{2,2}(0,1) : u = u(1-\cdot), u(0)= u(1) =h \}.
\end{equation}
  Here $\kappa$ denotes the signed curvature, $\mathrm{d}\mathbf{s}$ is the arclength element and $|\{u > 0 \}|$ is shorthand for the Lebesgue measure of $\{x \in (0,1) : u(x)> 0 \}$. Moreover, $\theta \in (0,\infty)$ is a real parameter which determines how much positivity is penalized. This \textit{elastic Alt-Caffarelli problem} has been studied before in \cite[Section 6]{DKV1} and relates to the interaction of elastic bending of a wire and gravitational force of a point mass on it. What makes the study interesting from the  point of view of calculus of variations is that minimizers have to balance out two competing interests: The second summand prefers a large nonpositivity set, whereas the first summand prefers small bending. However, in order to reach the zero level and turn back from there, bending is required. Our variational problem differs from the classical Alt-Caffarelli problem in two ways: the bending energy we use is nonlinear and of higher order, making the analysis involved even in one dimension. 
  
   Some classical challenges of  Alt-Caffarelli problems persist also in our setting. For example, due to the fact that the Lebesgue measure term is not a Frechét differentiable functional, Euler-Lagrange analysis cannot be used at the zero level $\{u = 0 \}$, which therefore forms a \textit{free boundary} of the problem. Alt-Caffarelli problems for higher order bending energies have raised some attention in the past, see e.g. \cite{DKV1,DKV2,Miura16,HG1,HG2,NahonLamboley,MariusAMPA,MariusIFB,MariusCVPD}. 
  Also closely related is the \textit{obstacle/confinement problem} for elastic graphs which has recently received a lot of attention, cf. \cite{GrunauOkabe1,GrunauOkabe2,DMN18,Bartels,ObstacleP,ObstacleProbl,Miura21,KY2021,DDObstacle}. Even 
a higher-dimensional analogue, the \textit{Willmore confinement problem} has been studied, cf. \cite{RoegerRupp}. 

Given that
\begin{equation}\label{eq:Etheta}
   \mathcal{B}(u):=  \int_{\mathrm{graph(u)}} \kappa^2 \; \mathrm{d}s = \int_0^1 \left(  \frac{u''(x)}{(1+ u'(x)^2)^\frac{3}{2}} \right)^2 \sqrt{1+ u'(x)^2} \; \mathrm{d}x = \int_0^1 \frac{u''(x)^2}{(1+ u'(x)^2)^\frac{5}{2}} \; \mathrm{d}x,
\end{equation}
the energy $\mathcal{E}_\theta$ is noncoercive, so that even existence of minimizers is not clear. We will see that the parameter $\theta$ plays a decisive role in the existence theory. It turns out that a crucial threshold for existence is given by $\theta= c_0^2$, where 
$
    c_0 = \int_{\mathbb{R}} \frac{1}{(1+s^2)^\frac{5}{4}} \; \mathrm{d}s  (\approx 2.39628). 
$
Our precise existence result is 
\begin{theorem}[Existence and Nonexistence]\label{thm1} Let $\mathcal{E}_\theta$ be as in \eqref{eq:Etheta}. 
    \begin{itemize}
        \item[(i)] Suppose that $\theta < c_0^2$. Then for each $h> 0$ the energy $\mathcal{E}_\theta$ possesses a minimizer in $\mathcal{A}_{\rm sym}(h)$. 
        \item[(ii)] Suppose that $\theta > c_0^2$. Then there exists $h_* > 0$ such that for all $h \geq h_*$ the energy $\mathcal{E}_\theta$ does not attain its infimum in $\mathcal{A}_{\rm sym}(h)$. 
    \end{itemize}
\end{theorem}
The theorem omits the case $\theta = c_0^2$. This case seems quite intriguing and not yet well-understood.  
It should be added that even in the case $\theta \geq c_0^2$ minimizers do exist for sufficiently small values of $h \in (0,\infty)$. 
Hence, it is a natural question to determine the optimal choice of $h_*$. An argument by Dall'Acqua and Deckelnick along the lines of \cite[Lemma 4.3]{DDObstacle} already gives a non-optimal value for $h_*$, namely one can choose $h_* = \frac{8}{c_0^2} (\approx 1.393)$. We recall this argument in Section \ref{sec:NonexistenceNonoptimal}, it relies on direct energy estimates and the Cauchy-Schwarz inequality. Examining optimal shapes more explicitly, we will however see that these estimates are not sharp for minimizers and the value of $h_*$ can be improved. 
\begin{theorem}[Improved nonexistence threshold]\label{thm2} Let $\theta > c_0^2$. Then the conclusion of Theorem \ref{thm1} (ii) holds true with 
\begin{equation}
    h_* = \sup_{b,s > 0} \frac{1}{\sqrt{\theta+b}} \int_{\theta/(\theta+b)}^1 \frac{s^2y}{(1+s^2y^2)^\frac{5}{4}} \frac{1}{\sqrt{1-y}} \; \mathrm{d}y. 
\end{equation}
This value is strictly smaller than $\frac{1}{\sqrt{\theta}}\frac{8}{c_0^2}.$  
\end{theorem}
Our proof of the fact that the supremum falls below $\frac{1}{\sqrt{\theta}}\frac{8}{c_0^2}$ is purely analytic and does not rely on numerics. The improved bound also shows that $h_* \to 0$ as $\theta \to \infty$, so the critical height depends vitally on $\theta$. It is also noteworthy that already for $\theta= c_0^2$ the quantity $\frac{1}{\sqrt{\theta}}\frac{8}{c_0^2}$ 
falls below optimal existence thresholds found in previous works for obstacle problems with \textit{cone obstacles} \cite{Miura21,KY2021}. This shows that existence theory for the Alt-Caffarelli problem is more involved than for the conical obstacle problem. 

Obtaining the improved height threshold requires precise understanding of \textit{optimal shapes}, i.e. one needs to study the structure of minimizers. We will now discuss some considerations about this structure which we shall prove rigorously in this article. 

If a minimizer $u^* \in \mathcal{A}_{\rm sym}(h)$ stays above the zero level, i.e $u^* > 0$ on $(0,1)$, then the Lebesgue measure term yields the energy value $\theta$ and thus $\mathcal{E}_\theta(u^*) \geq \theta$. Equality is actually possible if we look at the admissible function $\tilde{u} \equiv h \in \mathcal{A}_{\rm sym}(h)$. One infers that in this case $\mathcal{E}_\theta(u^*) = \theta$ and $u^* \equiv h$. For this reason the shape of minimizers is only interesting if $u^* \not> 0$. We speak of \textit{sign-changing minimizers}.  

\begin{figure}[h!]
    \centering
    \includegraphics[width=0.25\linewidth]{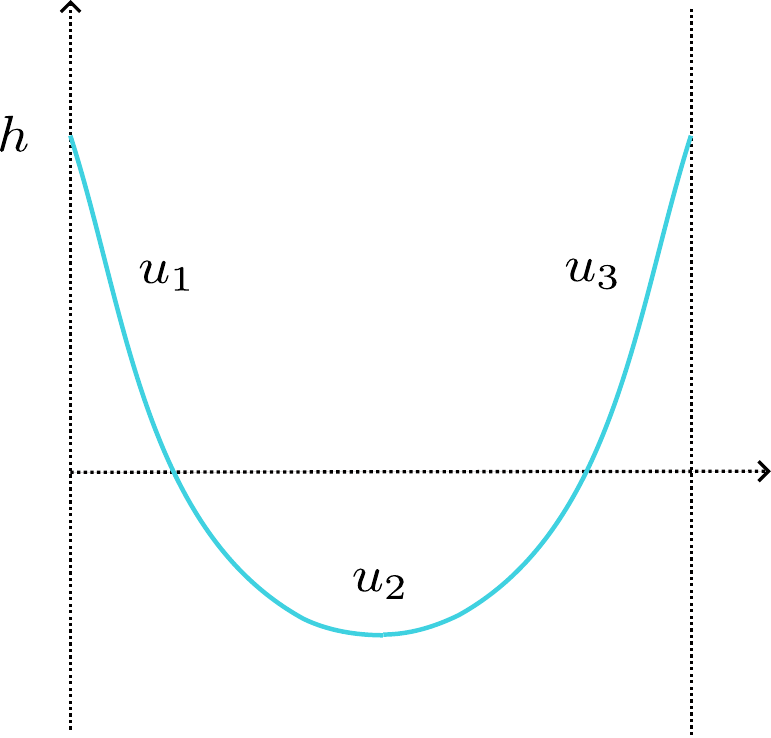}
    \caption{A schematic picture of a sign-changing minimizer $u^* = u_1 \oplus u_2 \oplus u_3$.}
    \label{fig:placeholder1}
\end{figure}

In this situation it seems reasonable to expect that minimizers will be composed of three graphs $u_1,u_2,u_3$, where $u_1,u_3>0$ connect the boundary level $h$ to the zero level and $u_2\leq 0$ (Figure \ref{fig:placeholder1}).
Due to reflection symmetry one has $u_3 = u_1(1-\cdot)$. Note that this decomposition is only possible if $u_*$ makes only two sign changes\footnote{To get this statement precise one needs the nonstandard convention that the sign of a real number $x$ is $1$ if $x> 0$ and $-1$ if $x \leq 0$}. Therefore, we need to answer the following question first
\begin{itemize}
    \item[(Q1)] Given a sign-changing minimizer, how many sign changes does it make? 
\end{itemize}
Due to the lack of a maximum principle for higher order problems this is a nontrivial question, even in one dimension. We will answer it in Section \ref{sec.Freeboundary} showing that only two sign-changes occur, justifying the above decomposition. 

Both $u_1$ and $u_3$ are critical points of $\mathcal{B}$, so called \textit{graphical elasticae}. The remaining function $u_2$ is also a graphical elastica on $\{ u_2 < 0 \}$, but a priori we do not know how large the set $\{ u_2 < 0 \}$ is. It is also possible that $u_2 \equiv 0$, or $u_2$ is composed of graphical elastica pieces with negative values and flat parts on the zero level. The third option seems unlikely for minimizers, but it seems unclear which of the first two versions is optimal. Hence we ask
\begin{itemize}
    \item[(Q2)] Is the middle part $u_2$ a flat line ($u_2 \equiv 0$)  or is $u_2$ a graphical elastica satisfying $u_2< 0$? 
\end{itemize}
\begin{figure}[h!]
    \centering
    \includegraphics[width=0.25\linewidth]{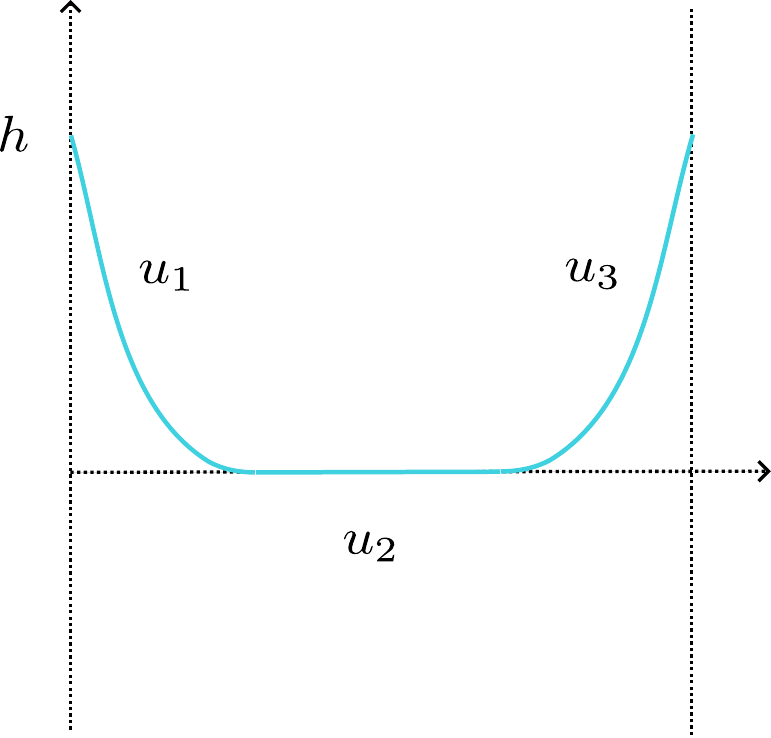}
    \caption{A competitor with a flat middle part, unlike in Figure \ref{fig:placeholder1}.}
    \label{fig:placeholder2}
\end{figure}

 At first sight it might seem convincing that $u_2$ is a line, as lines are absolute minimizers of $\mathcal{B}$. On the other hand, $u_2 \equiv 0$ would imply that $u_1$ and $u_3$ must approach the zero level in a flat way (Figure  \ref{fig:placeholder2}), which possibly requires extra bending. It turns out that this extra bending would cost too much energy and a non-flat configuration as e.g. in Figure \ref{fig:placeholder1} is preferred. 
 
This flatness question generally arises in higher order Alt-Caffarelli problems and can also be phrased in a more general way as follows. For a minimizer $u$, is the free boundary $\{u = 0\}$  reached in a flat way or is there slope at the free boundary?  This has been investigated in previous works by studying values of $|\nabla u|$ on $\{u = 0 \}$. In a linearized version of our functional, the (one-sided) biharmonic Alt-Caffarelli problem \cite{MariusAMPA,MariusCVPD}, it turns out that in dimension $2$ one has $|\nabla u| \neq 0$ on $\{ u= 0 \}$, i.e. one has slope everywhere at the boundary. The proof of this fact in \cite{MariusAMPA,MariusCVPD} works exclusively for dimension $2$ since the precise scaling of the Lebesgue measure is used in a decisive way. In this article we present a method that allows for the same conclusion in dimension $1$ (cf. Proposition \ref{prop:Gradientnonzero}). It relies on the interplay of the geometric scaling behavior of $\mathcal{B}$ and the observation that minimizers must be convex. Hence we can also answer (Q2) with the insight that the middle part is not flat. This in particular improves upon the analysis in \cite[Section 6, (6.4) and below]{DKV1}, where non-flatness is assumed but not verified.

Similar flatness questions also play a  role in the field of \textit{elastic networks} \cite{Network1,Network2, Network3}. Even in the most elementary case of \textit{theta networks of three elasticae} it is not understood whether minimizing networks contain a line. 
One major difference between our three-elastica system $u_1,u_2,u_3$ and an elastic theta network is however the \textit{angle condition}.  While we require that the tangents of the graphs match at their boundary, the network problem prescribes nonzero angles between the different curves (oftentimes an angle of $\frac{2\pi}{3}$).

Knowing that the middle part $u_2$ is not a line it seems reasonable to ask whether the remaining parts $u_1,u_3$ are lines:
\begin{itemize}
    \item[(Q3)] Are the graphs $u_1$ and $u_3$ affine-linear? 
\end{itemize}
What will however rule out this behavior is a regularity consideration: It will turn out that minimizers are $C^2$-regular. Therefore, affine linearity of $u_1$ (and thus of $u_3$) would imply that the curvature of $u_2$ must vanish at the boundary. We will see that this can not be the case, by having a closer look into symmetric graphical elasticae and therefore into the candidates for $u_2$ (Proposition \ref{prop:cshapeparameter}). Finally, we will even derive a stronger result, namely that each minimizer $u$ must be \textit{strictly convex} in $(0,1)$, ruling out line parts in all three sections at once. In particular, we can give a negative answer to (Q3).

The answers to the three previous questions (Q1),(Q2),(Q3) are summarized by the following \textit{structure theorem} for minimizers. 
 \begin{theorem}[Structure of minimizers]\label{thm3} 
     Suppose that $u \in \mathcal{A}_{\rm sym}(h)$ is a minimizer of $\mathcal{E}_\theta$ and $u \not \equiv h$. Then $u\in W^{3,\infty}(0,1)$ and $u''(0) = u''(1) =0$. Furthermore, $u$ is strictly convex on $(0,1)$ and there exists $a \in (0,\frac{1}{2})$ such that $\{ u = 0 \} = \{a,1-a\}$. Moreover, $u'(a), u'(1-a) \neq 0$ and $u\vert_{[0,a]},u \vert_{[a,1-a]}, u\vert_{[1-a,1]}$ are smooth graphical elasticae, none of which are lines.  
 \end{theorem}
Note that the curve in Figure \ref{fig:placeholder1} has all the properties in the above theorem and is therefore a reasonable candidate for a minimizer.
We will also see that the asserted $W^{3,\infty}$-regularity is optimal in the sense that minimizers can not lie in $W^{4,1}(0,1)$ unless they are constant (cf. Remark \ref{rem:nitchW41}). This regularity coincides with the obstacle problem of elastic graphs and might therefore not be surprising. 

Finally, we briefly comment on the symmetry assumption $u= u(1-\cdot)$. Due to the higher order nature and nonlinearity of the functional, symmetry of minimizers is nontrivial. More precisely, minimization in $\mathcal{A}(h):= \{ u \in W^{2,2}(0,1) : u(0)=u(1) = h \}$ could a priori yield different results. For small values of $h$ or small values of $\theta$ one can however obtain existence and symmetry of minimizers in $\mathcal{A}(h)$ along the lines of \cite[Section 6.1]{MariusAMOP}, using a nonlinear version of Talenti's symmetrization argument. We omit the argument here as it would essentially be a copy of \cite[Section 6.1]{MariusAMOP}. For larger values of $\theta$ and $h$ symmetry questions are to the best of the author's knowledge entirely open.

The article is structured as follows. In Section 2 we recall some known results, mainly from \cite{DDObstacle}. This allows us to prove Theorem \ref{thm1} with the nonoptimal value $h_* = \frac{8}{c_0^2}$. Section 3 is devoted to the description and properties of minimizers, i.e. the proof of Theorem  \ref{thm3}. Based on these findings we can improve the nonexistence threshold $h_*$ in Section 4 and prove Theorem \ref{thm2}.




\section{Preliminaries and previous work}
\subsection{Auxiliary functions}
In the following we recall a useful auxiliary function introduced in \cite{DDObstacle}.  Define the function $G: \mathbb{R} \to \mathbb{R}$ via 
\begin{equation}
    G(z) := \int_0^z \frac{1}{(1+s^2)^\frac{5}{4}} \; \mathrm{d}s. 
\end{equation}
Further we define 
\begin{equation}
    c_0 :=  2 \lim_{z \to \infty} G(z)   = \int_{\mathbb{R}} \frac{1}{(1+s^2)^\frac{5}{4}} \; \mathrm{d}s. 
\end{equation}
Note that $c_0 > 2$. Indeed, we have 
\begin{equation}\label{eq:c0istgroesser2}
  c_0 =  \int_{\mathbb{R}} \frac{1}{(1+s^2)^\frac{5}{4}} \; \mathrm{d}s > \int_{\mathbb{R}} \frac{1}{(1+s^2)^\frac{3}{2}} \; \mathrm{d}s  = \left[ \frac{s}{\sqrt{1+s^2}} \right]^\infty_{-\infty} = 2. 
\end{equation}
We will also show in  Lemma \ref{lem:cnullsquaresmaller6} that $c_0 < \sqrt{6}.$ Due to the fact that $G : \mathbb{R} \to (-\frac{c_0}{2},\frac{c_0}{2})$ is smooth and surjective and $G' > 0$ we have a smooth inverse function $G^{-1} :(-\frac{c_0}{2},\frac{c_0}{2}) \to \mathbb{R}$. For $c \in (0,c_0)$ we define
\begin{equation}\label{eq:uCAUx}
    u_c : [0,1] \to \mathbb{R} \quad u_c(x):= \frac{2}{c \sqrt[4]{1 + G^{-1}(\frac{c}{2}-cx)^2}} - \frac{2}{c \sqrt[4]{1 + G^{-1}(\frac{c}{2})^2}}.
\end{equation}
One readily checks that $u_c(0)= u_c(1)= 0, u_c( 1- \cdot) = u_c$ as well as $u_c \geq 0$ on $[0,1]$ and
\begin{align}
    u_c'(x) & =   \frac{1}{(1+ G^{-1}(\frac{c}{2}-cx)^2)^\frac{5}{4}} G^{-1}( \tfrac{c}{2}-cx) (G^{-1})'(\tfrac{c}{2}-cx)  \\ & = \frac{1}{(1+ G^{-1}(\frac{c}{2}-cx)^2)^\frac{5}{4}} 
  G^{-1}( \tfrac{c}{2}-cx)  \frac{1}{G'(G^{-1}(\tfrac{c}{2}-cx)^2)} = G^{-1}( \tfrac{c}{2}-cx).
\end{align}
This implies $(G \circ u_c')(x) = \frac{c}{2}-cx$ and hence
\begin{equation}\label{v_uisconstant}
    \frac{u_c''(x)}{(1+u_c'(x)^2)^\frac{5}{4}} = \frac{d}{dx} (G\circ u_c')(x)=  - c.
\end{equation}
The function $u_c$ defined as in \eqref{eq:uCAUx} will be a very useful auxiliary function throughout this article. 
\subsection{Basic existence theory}
In this section we  show that minimizers exist for $\theta < c_0^2$ and even for $\theta > c_0^2$ if $h$ is sufficiently small. Most of the techniques applied in this section are adaptations of \cite{DDObstacle} or \cite{DKV1} to our situation.
We first discuss an energy bound that is largely along the lines of \cite[Lemma 2.4]{DDObstacle}
\begin{prop}\label{prop:theta}
Let $\theta \in (0,\infty)$. Then
   $$ \inf_{\mathcal{A}_{\rm sym}(h)}  \mathcal{E}_\theta  \leq \min\{\theta, c_0^2\}.$$  
\end{prop}
\begin{proof}
Let $v : [0,1] \to \mathbb{R}$ be given by $v\equiv h$. Then clearly $v \in \mathcal{A}_{\rm sym}(h)$. Therefore, 
\begin{equation}\label{eq:energiekleinertheta}
  \inf_{\mathcal{A}_{\rm sym}(h)}  \mathcal{E}_\theta \leq   \mathcal{E}_\theta(v)= \int_0^1 0 \; \mathrm{d}x + \theta |(0,1)| = \theta. 
\end{equation}
Now let $\delta  \in (0, \frac{1}{2})$ and define (similar as in \cite[Lemma 2.4]{DDObstacle}), $v_\delta: [0,1] \to \mathbb{R}$ via 
\begin{equation}
    v_\delta(x) = \begin{cases}
        h- G^{-1}(\frac{c}{2}) x & x \in [0,\delta] \\
        h- G^{-1}(\frac{c}{2}) \delta  - (1-2\delta) u_c(\tfrac{x-\delta}{1-2\delta}) & x \in(\delta,1-\delta) \\
        h - G^{-1}(\frac{c}{2})(1-x) &x \in [1-\delta,1]
     \end{cases}.
\end{equation}
One readily checks that $v_\delta \in W^{2,2}(0,1)$, $v_\delta= v_\delta(1-\cdot)$ and $v_\delta(0)= v_\delta(1) = h$. Due to the fact that $\lim\limits_{c\uparrow c_0}G^{-1}(\frac{c}{2}) = \infty$ we can choose $\varepsilon > 0$ such that for all $c \in (c_0-\varepsilon,c_0)$ one has $\frac{2h}{G^{-1}(\frac{c}{2})}< \frac{1}{2}$. For $c\in (c_0-\varepsilon,c_0)$ choose $\delta_c := \frac{2h}{{G^{-1}(\frac{c}{2})}}$ and note that for $x \in(\delta_c,1-\delta_c)$
\begin{equation}
    v_{\delta_c}(x) \leq h - G^{-1}(\tfrac{c}{2}) \delta_c = - h < 0. 
\end{equation}
As a result $\{ v_{\delta_c} > 0 \} \subseteq (0,\delta_c] \cup [1-\delta_c,1)$. We infer for 
\begin{align}
    \mathcal{E}_\theta(v_{\delta_c}) & \leq  \int_0^1 \frac{v_{\delta_c}''(x)^2}{(1+ v_{\delta_c}'(x)^2)^{\frac{5}{4}}} \; \mathrm{d}x +\theta |\{ v_{\delta_c} > 0 \} | \leq   \int_{\delta_c}^{1-\delta_c} \frac{\frac{1}{(1-2\delta_c)^2} u_c''(\frac{x-\delta_c}{1-2\delta_c})^2}{( 1+ u_c'(\frac{x-\delta_c}{1-2\delta_c})^2)^\frac{5}{4}} \; \mathrm{d}x    + 2\theta\delta_c
     \\ & \leq \frac{1}{(1-2\delta_c)^2} c^2 (1-2\delta_c) + 2 \theta \delta_c = \frac{c^2}{1- \frac{4h}{G^{-1}(\frac{c}{2})}} + 4\theta h \frac{1}{G^{-1}(\frac{c}{2})}.
\end{align}
Letting $c \to c_0$ we have 
\begin{equation}
    \inf_{\mathcal{A}_{\rm sym}(h)} \mathcal{E}_\theta \leq \underset{c> c_0-\varepsilon}{\liminf_{c \to c_0} } \mathcal{E}_\theta(v_{\delta_c}) = c_0^2.
\end{equation}
The claim follows from this and \eqref{eq:energiekleinertheta}.
\end{proof}

If $h$ is small, the energy infimum can be a lot smaller than $\min\{ \theta, c_0^2 \}$ as the next proposition shows. 

\begin{prop}[Energy Bounds for small height]\label{prop:smallheight} Let $h < 1$. Then
    \begin{equation}
        \inf_{\mathcal{A}_{\rm sym}(h)} \mathcal{E}_\theta \leq (\theta + \tfrac{16}{\sqrt{5}})h^\frac{1}{3}.
    \end{equation}
\end{prop}
\begin{proof}
    Let $a := 2 h^\frac{1}{3}$. Notice that $\frac{a^2}{4} = h^\frac{2}{3} > h$.  Define 
       $u: (0,1) \to \mathbb{R},x \mapsto h-  \tfrac{a^2}{4} + a^2(x-\tfrac{1}{2})^2$. One readily checks that $u \in \mathcal{A}_{\rm sym}(h)$ and 
    $
       \left( u(x) > 0 \iff |x- \tfrac{1}{2}| > \sqrt{\tfrac{1}{4} - \tfrac{h}{a^2} } \right). 
   $
    Hence we have 
    \begin{equation}
        \mathcal{E}_\theta(u) =  \int_0^1 \frac{4a^2}{(1+ 4a^2(x-\frac{1}{2})^2)^\frac{5}{2}} \; \mathrm{d}x + 2 \theta \left( \tfrac{1}{2} - \sqrt{\tfrac{1}{4} - \tfrac{h}{a^2} } \right).
    \end{equation}
    Now we estimate $\frac{1}{2} - \sqrt{\frac{1}{4} - \frac{h}{a^2}} =  \frac{1}{\frac{1}{2} + \sqrt{\frac{1}{4} - \frac{h}{a^2}}}\frac{h}{a^2} \leq \frac{2h}{a^2}$ and, using that $a< 2,$ 
    \begin{equation}
        \int_0^1 \frac{4a^2}{(1+ 4a^2(x-\frac{1}{2})^2)^\frac{5}{2}} \; \mathrm{d}x =  2a \int_{-a}^{a}  \frac{1}{(1+z^2)^\frac{5}{2}} \; \mathrm{d}z \leq 2a \int_{-2}^2 \frac{1}{(1+z^2)^\frac{3}{2}} \; \mathrm{d}z = \frac{8}{\sqrt{5}}a. 
    \end{equation}
    Hence we obtain 
    $
        \mathcal{E}_\theta(u) \leq  4 \theta \frac{h}{a^2} + \frac{8}{\sqrt{5}}a  = (\theta + \tfrac{16}{\sqrt{5}})h^\frac{1}{3}.
   $
\end{proof}
The next existence criterion is a symbiosis of \cite[Lemma 4.2 and Lemma 2.5]{DDObstacle} and \cite[Lemma 2.1]{DKV1}.
\begin{prop}[Existence Criterion]\label{prop:excirterion}
Suppose that $\theta,h \in (0,\infty)$  are such that 
$
    \inf\limits_{\mathcal{A}_{\rm sym}(h)} \mathcal{E}_\theta< c_0^2.
$
Then  $\inf\limits_{\mathcal{A}_{\rm sym}(h)} \mathcal{E}_\theta$ is attained. 
\end{prop}
\begin{proof}
    Let $(u_n)_{n \in \mathbb{N}} \subseteq \mathcal{A}_{\rm sym}(h)$ be a minimizing sequence for $\mathcal{E}_\theta$. Since by assumption $\inf\limits_{\mathcal{A}_{\rm sym}(h)} \mathcal{E}_\theta < c_0^2$ we may assume that there exists $\delta  \in (0,c_0^2)$ such that 
    $\mathcal{E}_\theta(u_n) \leq c_0^2-\delta$ for all $n \in \mathbb{N}$. Note that due to the fact that $(u_n)_{n \in \mathbb{N}} \subseteq W^{2,2}((0,1))$ we have $(u_n')_{n \in \mathbb{N}} \subseteq C^0([0,1])$. 
    \begin{flushleft}
        \textbf{Claim 1.} $(u_n')_{n \in \mathbb{N}}$ is bounded in $C^0([0,1]).$
    \end{flushleft}
    To this end fix some $x_0 \in [0,\frac{1}{2}]$. Then for $n \in \mathbb{N}$ we have 
    \begin{align}
        c_0^2-\delta  & \geq \mathcal{E}_\theta(u_n) = \int_0^1 \frac{|u_n''(x)|^2}{(1+u_n'(x)^2)^\frac{5}{2}} \; \mathrm{d}x \geq \left( \int_0^1 \frac{|u_n''(x)| }{(1+ u_n'(x)^2)^\frac{5}{4}} \; \mathrm{d}x \right)^2   \geq \left( \int_{x_0}^{1-x_0} \frac{|u_n''(x)| }{(1+ u_n'(x)^2)^\frac{5}{4}} \; \mathrm{d}x \right)^2  
        \\ & \geq  \left(  \int_{x_0}^{1-x_0} |(G \circ u_n)'(x)| \; \mathrm{d}x \right)^2  \geq   \left\vert  \int_{x_0}^{1-x_0} (G \circ u_n)'(x) \; \mathrm{d}x \right\vert^2  =  [G(u_n'(x_0))-G(u_n'(1-x_0))]^2  \\ &=  [G(u_n'(x_0))+ G(-u_n'(1-x_0))]^2   = 4 G(u_n'(x_0))^2 = 4 G (|u_n'(x_0)|)^2, \label{eq:estiiderivative} 
    \end{align}
    where we have used in the last steps that $G$ is odd and  $u_n' = -u_n'(1-\cdot)$ due to the symmetry of $u_n$.
    We infer that 
    \begin{equation}
        |u_n'(x_0)| \leq G^{-1}\left( \tfrac{\sqrt{c_0^2-\delta}}{2} \right).
    \end{equation}
    Due to the symmetry of $u_n$ we find 
    $
        \sup_{[0,1]} |u_n'| = \sup_{[0,\frac{1}{2}]} |u_n'|  \leq  G^{-1}\left( \frac{\sqrt{c_0^2-\delta}}{2} \right).
    $
    The claim follows. 
    \begin{flushleft}
        \textbf{Claim 2.} $(u_n)_{n \in \mathbb{N}}$ is bounded in $W^{2,2}(0,1)$.
    \end{flushleft}
    Let $n \in \mathbb{N}$. For $x \in [0,1]$ we have 
    \begin{equation}
        u_n(x) = u_n(0) + \int_0^x u_n'(s) \; \mathrm{d}s  \leq h + \lVert u_n'\rVert_{C^0},
    \end{equation}
    so that we can estimate
   $
        \lVert u_n \rVert_{L^2} \leq \lVert u_n \rVert_{C^0} \leq h + \lVert u_n'\rVert_{C^0}
 $
    for all $n \in \mathbb{N}$. In particular $(u_n)_{n \in \mathbb{N}}$ is bounded in $L^2(0,1)$. Moreover, $\lVert u_n' \rVert_{L^2} \leq \lVert u_n' \rVert_{C^0}$, so that by Claim 1 $(u_n)_{n \in \mathbb{N}}$ is bounded in $W^{1,2}(0,1)$.
    Moreover, we compute with Proposition \ref{prop:theta}
    \begin{align}
        \lVert u_n'' \rVert_{L^2}^2  & = \int_0^1 u_n''(x)^2 \; \mathrm{d}x \leq (1+ \lVert u_n'\rVert_{C^0}^2)^\frac{5}{2}  \int_0^1 \frac{u_n''(x)^2}{(1+u_n'(x)^2)^\frac{5}{2}} \; \mathrm{d}x \\ & \leq (1+ \lVert u_n'\rVert_{C^0}^2)^\frac{5}{2} \mathcal{E}_\theta(u_n)  \leq  (1+ \lVert u_n'\rVert_{C^0}^2)^\frac{5}{2} (c_0^2-\delta),
    \end{align}
    which implies that $(u_n)_{n \in \mathbb{N}}$ is also bounded in $W^{2,2}(0,1)$. Hence there exists some $u_* \in W^{2,2}(0,1)$ such that after choosing a subsequence (no relabeling) $u_n \rightharpoonup u_*$ in $W^{2,2}(0,1)$. Due to the compact embedding $W^{2,2}(0,1) \hookrightarrow C^0([0,1])$ one readily checks that $u_* \in \mathcal{A}_{\rm sym}(h)$. 
    \begin{flushleft}
        \textbf{Claim 3.} $u_*$ minimizes $\mathcal{E}_\theta$ in $\mathcal{A}_{\rm sym}(h)$. 
    \end{flushleft}
    
   We first show that $|\{ u_* > 0 \}| \leq \liminf\limits_{n \to \infty} | \{ u_n > 0  \}|$. To this end let $x_0 \in (0,1)$ be such that $u_*(x_0) > 0$. Due to the fact that $u_n \to u_*$ in $C^0([0,1])$ there must exist some $n_0 \in \mathbb{N}$ such that $u_n(x_0) > 0$ for all $n \geq n_0$. This shows that 
   $
       \chi_{\{ u_* > 0 \}} \leq \liminf\limits_{n \to \infty} \chi_{\{ u_n  > 0 \}},
   $ where $\chi_A$ denotes the characteristic function of a set $A$. 
   This and Fatou's lemma yield 
   \begin{equation}
       |\{ u_* > 0 \} | = \int_0^1 \chi_{\{ u_* > 0 \}} \; \mathrm{d}x \leq  \int_0^1\liminf_{n \to \infty} \chi_{\{ u_n  > 0 \}} \leq  \liminf_{n \to \infty} \int_0^1 \chi_{\{ u_n  > 0 \}} = \liminf_{n \to \infty} |\{ u_n  > 0 \}| .
   \end{equation}
   Next we show that $\mathcal{B}(u_*) \leq \liminf\limits_{n \to \infty} \mathcal{B}(u_n)$. To this end we observe that due to the fact that $u_n'' \rightharpoonup u_*''$ in $L^2(0,1)$ we have 
   \begin{align}
       \mathcal{B}(u_*) & = \int_0^1 \frac{u_*''(x) u_*''(x)}{(1+u_*'(x)^2)^\frac{5}{2}} \; \mathrm{d}x = \lim_{n \to \infty} \int_0^1 \frac{u_*''(x) u_n''(x)}{(1+u_*'(x)^2)^\frac{5}{2}} \; \mathrm{d}x
       \\ &=  \lim_{n \to \infty} \int_0^1 \frac{u_*''(x) u_n''(x)}{(1+u_*'(x)^2)^\frac{5}{4} (1+ u_n'(x)^2)^\frac{5}{4}} \; \mathrm{d}x - \int_0^1 \frac{u_*''(x) u_n''(x)}{(1+u_*'(x)^2)^\frac{5}{4}} \left( \frac{1}{(1+u_*'(x)^2)^\frac{5}{4}} - \frac{1}{(1+u_n'(x)^2)^\frac{5}{4}} \right)  \; \mathrm{d}x. 
   \end{align}
   Using that $f: \mathbb{R} \to \mathbb{R}, z \mapsto \frac{1}{(1+z^2)^\frac{5}{4}}$ satisfies $|f'(z)| = \frac{5}{2} \frac{|z|}{(1+z^2)^\frac{9}{4}} \leq\frac{5}{4} \frac{1+z^2}{(1+z^2)^\frac{9}{4}} \leq \frac{5}{4}$ for all $z \in \mathbb{R}$ we may estimate the absolute value of the second summand from above  by 
   \begin{equation}
        \lVert u_*'' \rVert_{L^2} \lVert u_n'' \rVert_{L^2} \frac{5}{4} \lVert u_n' - u_*'\rVert_{C^0} \xrightarrow{n \to \infty} 0,
   \end{equation}
   due to the fact that the embedding $W^{2,2}(0,1) \hookrightarrow C^1([0,1])$ is compact. Using this and the Cauchy-Schwarz inequality in the first summand we obtain 
   \begin{equation}
       \mathcal{B}(u_*) = \liminf_{n \to \infty} \sqrt{\int_0^1 \frac{u_*''(x)^2}{(1+ u_*'(x)^2)^\frac{5}{2}} \; \mathrm{d}x}  \sqrt{\int_0^1 \frac{u_n''(x)^2}{(1+ u_n'(x)^2)^\frac{5}{2}} \; \mathrm{d}x}  =   \liminf_{n \to \infty} \sqrt{\mathcal{B}(u_*)} \sqrt{\mathcal{B}(u_n)}.
   \end{equation}
   The lower semicontinuity of $\mathcal{B}$ at $u_*$ follows immediately. Finally, we have $u_* \in \mathcal{A}_{\rm sym}(h)$ and
    \begin{align}
       \mathcal{E}_\theta(u_*)  & = \mathcal{B}(u_*) + \theta|\{ u_* > 0 \}| \leq \liminf_{n \to \infty} \mathcal{B}(u_n) +  \theta\liminf_{n \to \infty} |\{ u_n > 0 \}|  \\ & \leq \liminf_{n \to \infty} ( \mathcal{B}(u_n) +\theta |\{ u_n > 0 \}|) = \liminf_{n \to \infty} \mathcal{E}_\theta(u_n) = \inf_{\mathcal{A}_{\rm sym}(h)} \mathcal{E}_\theta. \qedhere
   \end{align} 
\end{proof}

\begin{cor}\label{cor:theta<c_0^2}
    Suppose that $\theta< c_0^2$. Then $\inf\limits_{\mathcal{A}_{\rm sym}(h)} \mathcal{E}_\theta$ is attained for any $h> 0$.
\end{cor}
\begin{proof}
If $\theta < c_0^2$ then Proposition \ref{prop:theta} yields 
$
    \inf\limits_{\mathcal{A}_{\rm sym}(h)} \mathcal{E}_\theta \leq \theta < c_0^2,
$
so that Proposition \ref{prop:excirterion} yields the result. 
\end{proof}
\begin{remark}
    If $\theta < c_0^2$ and $h>0$ is suitably large then it is expected that $u \equiv h$ is the unique minimizer of the problem. Indeed, proceeding as in \eqref{eq:estiiderivative} we can estimate for a minimizer $u \in \mathcal{A}_{\rm sym}(h)$ that $\mathcal{E}_\theta(u) \geq 4 G(\sup_{[0,1]} |u'|)^2$. Estimating $\mathcal{E}_\theta(u) \leq \theta$ by Proposition \ref{prop:theta} we obtain $\sup_{[0,1]} |u'| \leq G^{-1}(\frac{\sqrt{\theta}}{2})$. Now, if $u \not \equiv h$ then $u$ must change sign and thus $\sup_{[0,1]} |u'| \geq h$, yielding $h \leq G^{-1}(\frac{\sqrt{\theta}}{2})$. The contrapositive of this statement is that $h > G^{-1}(\frac{\sqrt{\theta}}{2})$ implies that minimizers are constant. 
\end{remark}
\begin{remark}
    From Proposition \ref{prop:excirterion} and Proposition \ref{prop:smallheight} follows also that regardless of the value of $\theta$ we have existence of minimizers if $h$ is small enough, more precisely $h < \min\{1, \frac{c_0^6}{(\theta + \frac{16}{\sqrt{5}})^3} \}$.  
\end{remark}

\subsection{Nonexistence}\label{sec:NonexistenceNonoptimal}
Next we investigate criteria for nonexistence of minimizers. By Corollary \ref{cor:theta<c_0^2} this is only possible if $\theta \geq c_0^2$. We proceed along the lines of \cite[Lemma 4.3]{DDObstacle} to derive 
\begin{prop}\label{prop:NaiveNonexistence}
    Suppose that $\theta > c_0^2$. Then for each $h > \frac{8}{c_0^2}$ the infimum $\inf\limits_{\mathcal{A}_{\rm sym}(h)} \mathcal{E}_\theta$ is not attained.
\end{prop}
\begin{proof}
    Let $h > \frac{8}{c_0^2}$. Assume the opposite, i.e. there exists a minimizer $u \in \mathcal{A}_{\rm sym}(h)$. Note that $u$ must attain nonpositive values, since $u> 0$ on $(0,1)$ would imply $\mathcal{E}_\theta(u) = \mathcal{B}(u) + \theta|(0,1)| \geq \theta > c_0^2$, which contradicts the minimality of $u$ due to Proposition \ref{prop:theta}. Hence, due to the intermediate value theorem and symmetry of $u$, there must exist $x_0 \leq \frac{1}{2}$ such that $u(x_0) =0$. We calculate now for $t \in [0,\frac{1}{2})$ 
    \begin{align}
        |G(u'(t))|  & = |G(u'(t)) - G(u'(\tfrac{1}{2}))| = \left\vert \int_{\frac{1}{2}}^t \frac{u''(s)}{(1+u'(s)^2)^\frac{5}{4}} \; \mathrm{d}s  \right\vert \\& \leq \int_{t}^\frac{1}{2} \frac{|u''(s)|}{(1+u'(s)^2)^\frac{5}{4}} \; \mathrm{d}s  \leq \sqrt{\tfrac{1}{2}-t} \sqrt{ \int_{t}^\frac{1}{2} \frac{u''(s)^2}{(1+u'(s)^2)^\frac{5}{2}} \; \mathrm{d}s}
        = \sqrt{\tfrac{1}{2}-t} \sqrt{\tfrac{1}{2} \mathcal{B}(u) }.
    \end{align}
    Estimating $\mathcal{B}(u) \leq \mathcal{E}_\theta(u) \leq c_0^2$ with the aid of Proposition \ref{prop:theta} we obtain $G(|u'(t)|) = |G(u'(t))| \leq \frac{c_0}{2} \sqrt{1-2t}$ for all $t \in (0,\frac{1}{2})$. As a consequence we find $|u'(t)| \leq G^{-1}( \frac{c_0}{2} \sqrt{1-2t})$ for all $t \in (0, \frac{1}{2})$. Thus we can estimate, following the lines of \cite[end of proof of Lemma 4.3]{DDObstacle}
    \begin{align}
        h & = u(0) - u(x_0) = \int_{x_0}^0 u'(t) \; \mathrm{d}t \leq \int_0^{x_0} |u'(t)| \; \mathrm{d}t \leq  \int_0^{\frac{1}{2}} G^{-1}( \tfrac{c_0}{2} \sqrt{1-2t}) \; \mathrm{d}t
        \\ & = \frac{4}{c_0^2} \int_0^\frac{c_0}{2} \tau G^{-1}(\tau) \; \mathrm{d}\tau = \frac{4}{c_0^2} \int_0^\infty G(s) \frac{s}{(1+s^2)^\frac{5}{4}} \; \mathrm{d}s =- \frac{8}{c_0^2} \int_0^\infty G(s) \frac{d}{ds} \left( \frac{1}{(1+s^2)^\frac{1}{4}} \right) \; \mathrm{d}s 
        \\ & = \frac{8}{c_0^2} \int_0^\infty G'(s) \frac{1}{(1+s^2)^\frac{1}{4}} \; \mathrm{d}s = \frac{8}{c_0^2} \int_0^\infty \frac{1}{(1+s^2)^\frac{3}{2}} \; \mathrm{d}s = \frac{8}{c_0^2} \left[ \frac{s}{\sqrt{1+s^2}}\right]^{s=\infty}_{s=0} = \frac{8}{c_0^2}.
    \end{align}
    This yields $h \leq \frac{8}{c_0^2}$, a contradiction. 
\end{proof}

With this information at hand we can give a proof of Theorem \ref{thm1} with a nonoptimal value of $h_*$.
\begin{proof}[Proof of Theorem \ref{thm1}.]
    Part (i) follows from Corollary \ref{cor:theta<c_0^2} and Part (ii) can be obtained by choosing $h_* = \frac{8}{c_0^2}$ and using Proposition \ref{prop:NaiveNonexistence}.
\end{proof}

As announced in the introduction we will improve this height threshold $h_*$ later, once we have a more explicit understanding of the shape of minimizers.

\section{Description of minimizers}
This section is devoted to the study of properties of minimizers i.e. to the proof of Theorem \ref{thm3}. First, we will use variational methods to derive a \textit{measure-valued Euler Lagrange equation}. Afterwards, convexity considerations and the geometry of the problem will allow us to characterize the measure term in this Euler Lagrange equation explicitly up to a two-parameter family. This information is already enough to characterize the free boundary and derive the desired strict convexity.
\subsection{A measure-valued Euler-Lagrange equation}
Let $u \in W^{2,2}(0,1)$. We define $v_u \in L^2(0,1)$ via the formula
\begin{equation}
 v_u(x) := \frac{u''(x)}{(1+u'(x)^2)^\frac{5}{4}} \qquad \textrm{a.e. $x \in (0,1)$.}
\end{equation}
This definition allows us to describe critical points of $\mathcal{B}$, so called \textit{graphical elasticae}, in a compact way. 
\begin{lemma}\label{lem:33}
    Let $\varphi \in W^{2,2}(0,1)$. Then the Frechét derivative of $\mathcal{B} : W^{2,2}(0,1) \to \mathbb{R}, u \mapsto \mathcal{B}(u)$ satisfies 
    \begin{equation}
        D\mathcal{B}(u) (\varphi) = 2 \int_0^1  v_u(x) \frac{d}{dx} \frac{\varphi'(x)}{(1+ u'(x)^2)^\frac{5}{4}} \; \mathrm{d}x. 
    \end{equation}
\end{lemma}
\begin{proof}
    One readily computes that for $\varphi$ as in the statement one has
    \begin{align}
        D\mathcal{B}(u)(\varphi)   & =  \int_0^1  \left( \tfrac{ 2 u''(x) \varphi''(x)}{(1+u'(x)^2)^\frac{5}{2}} - \tfrac{5 u''(x)^2 u'(x)\varphi'(x) }{(1+u'(x)^2)^\frac{7}{2}} \right)  \; \mathrm{d}x 
      \\&  =  2 \int_0^1 \tfrac{u''(x)}{(1+u'(x)^2)^\frac{5}{4}}  \left( \tfrac{\varphi''(x)}{(1+u'(x)^2)^\frac{5}{4}} - \tfrac{5}{2}\tfrac{ u''(x) u'(x)\varphi'(x) }{(1+u'(x)^2)^\frac{9}{4}}  \right)   \; \mathrm{d}x  = 2 \int_0^1 v_u(x) \frac{d}{dx} \left(\tfrac{\varphi'(x)}{(1+ u'(x)^2)^\frac{5}{4}}\right) \; \mathrm{d}x. \qedhere
    \end{align}
\end{proof}

\begin{prop}\label{prop:preevious}
  Let $u \in \mathcal{A}_{\rm sym}(h)$ be a minimizer for $\mathcal{E}_\theta$. Further let $a,b \in (0,1), a<b$ be such that $u$ does not have a zero in $(a,b)$. 
  Then $ u \in C^\infty([a,b])$ and $v_u$  satisfies $\frac{d}{dx} \frac{v_u'}{(1+u'^2)^\frac{5}{4}} = 0$ on $[a,b]$. 
\end{prop}
\begin{proof}
Let $\eta \in C_0^\infty( a,b )$. Define $\tilde{\eta} := \frac{1}{2}( \eta + \eta ( 1- \cdot))$. Observe that $\tilde{\eta} \in  C_0^\infty((a,b)\cup (1-b,1-a))$ and $\lVert \tilde{\eta} \rVert_{C^0} \leq \lVert \eta \rVert_{C^0}$. 
Now $|u|$ attains a minimum on the compact set $\mathrm{spt}(\tilde{\eta})$ and since $|u| > 0$ on $(a,b) \cup (1-b,1-a)$ we have $ \delta := \min_{\mathrm{spt}(\tilde{\eta})} |u| > 0$. We next look at the variation $u_t := u+ t \tilde{\eta}$. One readily checks that for all $t \in (-\frac{\delta }{\lVert \tilde{\eta} \rVert_{C^0} }, \frac{\delta }{\lVert \tilde{\eta} \rVert_{C^0} }) =:I $ one has $\{ u + t \eta > 0 \} = \{ u > 0\}$. In particular, for $t \in I$ we have $\mathcal{E}_\theta(u+t\tilde\eta) = \mathcal{B}(u+t\tilde\eta) + \theta |\{ u > 0 \}|$. This expression is differentiable at $t = 0$. Minimality of $u$ yields
    \begin{equation}\label{eq:DB=0}
        0 = \frac{d}{dt} \mathcal{E}_\theta(u+t\tilde\eta) \Big\vert_{t = 0} = D\mathcal{B}(u)(\tilde\eta). 
    \end{equation}
Due to the symmetry of $u$ one readily checks that $D \mathcal{B}(u) (\tilde{\eta}) = D\mathcal{B}(u)(\eta)$ and hence 
\begin{equation}
    D\mathcal{B}(u)(\eta) = 0 \qquad \textrm{for all $\eta \in C_0^\infty((a,b))$.}
\end{equation}
    One can now argue as in \cite[Proposition 3.2(a)]{DDObstacle} to show that $u \in C^\infty([a,b])$. Given this and \eqref{eq:DB=0} we may use Lemma \ref{lem:33} to conclude that for every $\eta \in C_0^\infty((a,b))$ we have 
    \begin{equation}
        0 = D \mathcal{B}(u)(\eta) = 2 \int_a^b  v_u(x) \frac{d}{dx} \frac{\eta'(x)}{(1+ u'(x)^2)^\frac{5}{4}} \; \mathrm{d}x = 2 \int_a^b  \eta(x)  \frac{d}{dx}  \frac{v_u' (x)}{(1+ u'(x)^2)^\frac{5}{4}} \; \mathrm{d}x. 
    \end{equation}
    The fundamental lemma of calculus of variations yields $\frac{d}{dx} \frac{v_u'}{(1+u'^2)^\frac{5}{4}} = 0$ on $(a,b)$ (which holds also true on $[a,b]$ by continuity of the respective derivatives). 
\end{proof}

\begin{remark}\label{eq:Cinfregonuneq0}
    An important consequence of the previous proposition is that $u \in C^\infty(\{ u \neq 0 \})$ and $\frac{d}{dx} \frac{v_u'}{(1+u'^2)^\frac{5}{4}}= 0$ on $\{u \neq 0\}$. This is the so-called \textit{graphical elastica equation}.   
\end{remark}

\begin{lemma}\label{lem:36}
    There exists $\delta \in (0,\frac{1}{2})$ such that $u \in C^\infty([0,\delta] \cup [1-\delta,1])$ and $u''(0)=u''(1)=0$. 
\end{lemma}
\begin{proof}
By symmetry it suffices to show that for some $\delta \in (0,\frac{1}{2})$ we have $u \in C^\infty([0,\delta])$ and $u''(0)= 0$.  Due to continuity and the fact that $u(0)= h> 0$ there exists $\delta > 0$ such that $(0,\delta) \subset \{ u \neq 0\}$. In particular, the fact that $u \in C^\infty([0,\delta])$ follows from \eqref{prop:preevious}. Now let $\eta \in C_0^\infty([0,\delta))$ be such that $\eta(0)=0$ (but not necessarily $\eta'(0)=0$). Define $\tilde{\eta}:= \frac{1}{2}(\eta + \eta(1-\cdot))$. Then, for $|t| > 0$ small enough we have that $u_t := u + t \tilde\eta \in \mathcal{A}_{\rm sym}(h)$ as well as $|\{ u +t \tilde\eta > 0 \}| = |\{ u > 0\}|$ (arguing as in the proof of Proposition \ref{prop:preevious}). This yields again that for $|t|>0$ small enough one has $\mathcal{E}_\theta(u+ t \eta)= \mathcal{B}(u+t\eta) + \theta |\{ u > 0 \}|$. We infer differentiability at $t= 0$ and compute again due to minimality 
\begin{equation}
    0 = \frac{d}{dt} \mathcal{E}_\theta(u+ t \tilde{\eta}) \Big\vert_{t = 0} = D \mathcal{B}(u)(\tilde{\eta}).
\end{equation}
Symmetry of $u$ yields again $D \mathcal{B}(u)(\eta)= 0$. Hence Lemma \ref{lem:33} yields that 
\begin{align}
  0  & = \int_0^\delta v_u \frac{d}{dx} \frac{\eta'}{(1+u'^2)^\frac{5}{4}} \; \mathrm{d}x = \left[ \frac{v_u}{(1+u'^2)^\frac{5}{4}}  \eta' \right]_0^\delta - \int_0^\delta \frac{v_u'}{(1+u'^2)^\frac{5}{4}} \eta' \; \mathrm{d}x
  \\ & = \left[ \frac{v_u}{(1+u'^2)^\frac{5}{4}}  \eta' \right]_0^\delta - \left[ \frac{v_u'}{(1+u'^2)^\frac{5}{4}} \eta \right]_0^\delta  + \int_0^\delta \frac{d}{dx} \frac{v_u'}{(1+u'^2)^\frac{5}{4}} \eta \; \mathrm{d}x= \left[ \frac{v_u}{(1+u'^2)^\frac{5}{4}}  \eta' \right]_0^\delta,
\end{align}
 where we used in the last step that $\eta(0)= \eta(\delta) =0$ and that by Proposition \ref{prop:preevious} $\frac{d}{dx} \frac{v_u'}{(1+u'^2)^\frac{5}{4}} =0$ on $[0,\delta]$. Evaluating this expression we obtain $0 = \frac{v_u(0)}{(1+u'(0)^2)^\frac{5}{4}} = \frac{u''(0) }{(1+u'(0)^2)^\frac{5}{2}}$. This implies $u''(0)= 0$ as claimed. 
\end{proof}

\begin{lemma}\label{lem:RadonMeas}
    Let $u \in \mathcal{A}_{\rm sym}(h)$ be a minimizer of $\mathcal{E}_\theta$. Then there exists a finite Radon measure $\mu$ supported on $\{u = 0 \}$ such that for all $\varphi \in W^{2,2}(0,1) \cap W_0^{1,2}(0,1)$
    \begin{equation}\label{eq:todo}
        D \mathcal{B}(u)(\varphi) = - \int_{\{u = 0 \} } \varphi \; \mathrm{d}\mu . 
    \end{equation}
   In particular, $u \in W^{3,\infty}(0,1)$ and $v_u \in W^{1,\infty}(0,1)$.
\end{lemma}
\begin{proof}
    Let $\varphi \in C_0^\infty((0,1))$ be such that $\varphi \leq 0$. Define $\tilde{\varphi} := \frac{1}{2}(\varphi+ \varphi(1-\cdot))$
    Then for all $t \geq 0$ one has $|\{ u + t \tilde{\varphi}> 0 \}| \leq |\{ u > 0 \}|$. Due to the minimality of $u$ and the fact that $u+ t \tilde{\varphi} \in \mathcal{A}_{\rm sym}(h)$ 
    we have for any $t \geq 0$
    \begin{equation}
        \mathcal{B}(u) + \theta|\{ u > 0 \}| = \mathcal{E}_\theta(u) \leq \mathcal{E}_\theta(u+t\tilde{\varphi})=   \mathcal{B}(u+t \varphi) + \theta|\{ u+t \tilde{\varphi} > 0 \}| \leq \mathcal{B}(u+t \tilde{\varphi}) + \theta|\{ u > 0 \}|. 
    \end{equation}
    Subtracting $\theta |\{ u > 0 \}|$ from both sides we obtain $\mathcal{B}(u) \leq \mathcal{B}(u+t\tilde{\varphi})$. This implies 
    \begin{equation}
        D \mathcal{B}(u) (\tilde{\varphi}) = \lim_{t \to 0} \frac{\mathcal{B}(u+t\tilde{\varphi})- \mathcal{B}(u)}{t} = \underset{t >0 }{\lim_{t \to 0}} \frac{\mathcal{B}(u+t\tilde{\varphi})- \mathcal{B}(u)}{t}  \geq 0.
    \end{equation}
Again by symmetry of $u$ one readily checks that $D\mathcal{B}(u)(\tilde{\varphi})=D\mathcal{B}(u)(\varphi)$. Hence we infer 
\begin{equation}
    D \mathcal{B}(u) (\varphi) \geq 0 \quad \textrm{for all $\varphi \in C_0^\infty(0,1)$.}
\end{equation}
 Now define $L : C_0^\infty((0,1)) \to \mathbb{R}$ by $L(\psi):= D \mathcal{B}(u)(-\psi)$. The above computation shows that for each $\psi \in C_0^\infty((0,1))$ such that $\psi \geq 0$ we have $L(\psi) \geq 0$. By the Riesz-Markov-Kakutani theorem (cf. \cite[Theorem 1.39]{EVGARIEPY2015}) there exists a Radon measure $\mu$ such that 
    \begin{equation}
        L(\psi) = \int \psi \; \mathrm{d}\mu \quad \textrm{for all $\psi \in C_0^\infty(0,1)$.}
    \end{equation}
    In particular, 
    \begin{equation}\label{eq:eqphidmu}
        D \mathcal{B}(u)(\varphi) = - \int \varphi \; \mathrm{d}\mu. 
    \end{equation}
    \begin{flushleft}
        \textbf{Claim 1.} $\mu$ is finite and supported on $\{ u = 0 \}$.
    \end{flushleft}
    To this end let $\varphi \in C_0^\infty( \{ u \neq 0 \}).$ Then we can use \eqref{eq:eqphidmu} and Remark \ref{eq:Cinfregonuneq0} to find that 
    \begin{equation}
        -\int \varphi \; \mathrm{d}\mu = D \mathcal{B}(u)(\varphi) = 2 \int_0^1 v_u(x) \frac{d}{dx} \frac{\varphi'(x)}{(1+u'(x)^2)^\frac{5}{4}} \; \mathrm{d}x = 2 \int_0^1 \frac{d}{dx}\frac{v_u'(x)}{(1+u'(x)^2)^\frac{5}{4}} \varphi(x) \; \mathrm{d}x = 0 . 
    \end{equation}
    Due to the fact that $\{ u \neq 0 \}$ is open we can approximate for each compact set $K \subset \{ u \neq 0 \}$ the characteristic function $\chi_K$ uniformly by functions in $C_0^\infty(\{ u \neq 0 \})$. This and the fact that $\mu$ is a Radon measure yields $\mu(\{ u \neq 0 \})= 0$. The finiteness of $\mu$ follows from the fact that $\mu$ is a Radon measure and $\{u = 0 \} \subset (0,1)$ is compact (since $u(0)=u(1)= h>0$ and hence $u > 0$ in a neighborhood of $0$ and $1$).  
\begin{flushleft}
    \textbf{Claim 2.} \eqref{eq:todo} holds true for all $\varphi \in W^{2,2}(0,1) \cap W_0^{1,2}(0,1)$.
\end{flushleft}
Let $\varphi \in W^{2,2}(0,1) \cap W_0^{1,2}(0,1)$ be arbitrary. Further choose $\eta \in C_0^\infty(0,1)$ such that $\eta \equiv 1$ on an open neighborhood of $\{u = 0 \}$. Then $\eta \varphi \in W^{2,2}(0,1)$ has compact support. For this reason there exists $(\varphi_\varepsilon)_{\varepsilon > 0} \subset C_0^\infty(0,1)$ such that $\varphi_\varepsilon \rightarrow \eta \varphi$ in $W^{2,2}(0,1)$ and thus also uniformly on $[0,1]$ as $\varepsilon \to 0$. Using the finiteness of $\mu$ we have 
\begin{align}
    -\int \varphi \; \mathrm{d}\mu & = -\int_{\{ u = 0 \}} \varphi \; \mathrm{d}\mu = -\int_{\{ u = 0 \}} \varphi  \eta \; \mathrm{d}\mu = \lim_{\varepsilon \to 0} -\int_{\{ u = 0 \}} \varphi_\varepsilon \; \mathrm{d}\mu 
    \\ & = \lim_{\varepsilon \to 0} D \mathcal{B}(u) (\varphi_\varepsilon) =  D \mathcal{B}(u) ( \varphi \eta), \label{eq:restvomfest}
\end{align}
where we used continutity of $D \mathcal{B}(u): W^{2,2}(0,1) \to \mathbb{R}$ in the last step (cf. Lemma \ref{lem:33}). 
Now $\varphi (1-\eta)$ is supported on $[0,1] \setminus O$ where $O$ is an open neighborhood of $\{ u = 0 \}$ and vanishes on the endpoints $0,1$. This implies that there exists some $\delta > 0$ such that $|u| > \delta$ on $\mathrm{spt}(\varphi (1-\eta))$. As a consequence one infers, looking at the perturbation $u_t := u + t \varphi (1-\eta)$ and proceeding as in the proof of Proposition \ref{prop:preevious}, that 
\begin{equation}
    D \mathcal{B}(u)(\varphi(1-\eta))= 0.
\end{equation}
This, \eqref{eq:restvomfest}, and linearity of $D \mathcal{B}(u)$ yields that 
\begin{equation}
    D \mathcal{B}(u)(\varphi) =- \int \varphi \; \mathrm{d}\mu = - \int_{\{ u = 0 \}} \varphi \; \mathrm{d}\mu  .
\end{equation}
\begin{flushleft}
        \textbf{Claim 3.} $u \in W^{3,\infty}(0,1)$ (and as a consequence $v_u \in W^{1,\infty}(0,1)$)
    \end{flushleft}
   Once \eqref{eq:eqphidmu} is established, this follows precisely the along lines of \cite[Theorem 5.1]{DDObstacle}.
\end{proof}

One further identity that is based on the principle of inner variation will be of help --- especially for characterizing the Radon measure $\mu$. 

\begin{lemma}\label{eq:LemmaRadonMeasExplicit}
    The Radon measure $\mu$  in Lemma \ref{lem:RadonMeas} satisfies 
    \begin{equation}
        \int_{\{u = 0 \}} u' \eta \; \mathrm{d}\mu = - \theta \int_{\{ u > 0 \}} \eta'(x) \; \mathrm{d}x  \quad \textrm{for all $\eta \in C_0^\infty(0,1)$.}
    \end{equation}
\end{lemma}
\begin{proof}
    Let $\eta \in C_0^\infty(0,1)$. Observe that $\tilde{\eta} := \frac{1}{2}(\eta - \eta(1-\cdot))$ also lies in $C_0^\infty((0,1))$ (the minus sign is intended). Then for $|t|$ small enough the map $\phi_t: (0,1) \to (0,1), x \mapsto x+ t \tilde\eta(x)$ is a smooth diffeomorphism that maps $(0,1)$ onto itself and even equals the identity outside of a compact set $K \subset (0,1)$. Now define $u_t:= u \circ \phi_t$. One readily checks that $t \mapsto u_t$ lies in $C^1((-\varepsilon,\varepsilon),W^{2,2}(0,1))$ for some $\varepsilon > 0$. Notice carefully that this uses the $W^{3,\infty}$-regularity result  from Lemma \ref{lem:RadonMeas}. Indeed, to ensure that the time derivative $\partial_t u_t = (u' \circ \phi_t ) \cdot \partial_t \phi_t = (u' \circ \phi_t) \cdot \tilde\eta$ lies in $W^{2,2}(0,1)$ one needs that $u' \in W^{2,2}(0,1)$, which is true only due to Lemma \ref{lem:RadonMeas}.
Due to this fact $t \mapsto \mathcal{B}(u_t)$ is differentiable at $t= 0$ and 
\begin{equation}
    \frac{d}{dt} \mathcal{B}(u_t) \Big\vert_{t= 0} = D \mathcal{B}(u) (\partial_t u_t \vert_{t= 0}) = D \mathcal{B}(u)(u' \tilde{\eta}). 
\end{equation}
Due to symmetry of $u$ we have $D \mathcal{B}(u)(u' \eta) =D \mathcal{B}(u)( \frac{1}{2} (u' \eta+ (u' \eta)(1-\cdot) ) $. Given that $u'(1-\cdot)=-u'$ we obtain $D \mathcal{B}(u)(u' \eta)= D \mathcal{B}(u)(u' \tilde\eta)$. We infer
\begin{equation}\label{eq:deri1}
    \frac{d}{dt} \mathcal{B}(u_t) \Big\vert_{t= 0} = D \mathcal{B}(u)(u' \eta) =  - \int_{\{u = 0\}} u' \eta \; \mathrm{d}\mu,
\end{equation}
where we used Lemma \ref{lem:RadonMeas} in the last step. 
 Now for $|t|$ small enough we have 
    \begin{align}
        |\{ u_t > 0 \}| = |\phi_t^{-1} \{ u > 0 \} | = \int_{\{u > 0\}} |(\phi_t^{-1})'(x)| \; \mathrm{d}x = \int_{\{u > 0\}} \frac{1}{|\phi_t'(\phi_t^{-1}(x))|} \; \mathrm{d}x  = \int_{\{u > 0\}} \frac{1}{1+ t \tilde\eta'(\phi_t^{-1}(x))} \; \mathrm{d}x.
    \end{align}
   Given this, one readily checks that $t \mapsto \{u_t > 0 \}$ is differentiable at $t= 0$ and 
   \begin{equation}
     \frac{d}{dt}  |\{ u_t > 0 \}| = -\int_{\{ u > 0 \}} \tilde\eta'(x) \; \mathrm{d}x = - \frac{1}{2}\int_{\{u > 0\}} (\eta'(x) + \eta'(1-x)) = -\int_{\{ u > 0 \}} \eta'(x) \; \mathrm{d}x,
   \end{equation}
   where we used the symmetry of $\{ u > 0\}$ in the last step. 
   We have derived with this and \eqref{eq:deri1} that
   \begin{equation}\label{eq:ddtequt}
       \frac{d}{dt} \mathcal{E}_\theta(u_t) \Big\vert_{t = 0} = - \int_{\{ u = 0 \}} u' \eta \; \mathrm{d}\mu - \theta\int_{\{ u > 0 \} } \eta'(x) \; \mathrm{d}x. 
   \end{equation}
   Now notice that for all $|t|$ small enough we have $u_t \in \mathcal{A}_{\rm sym}(h)$. Indeed, the boundary values are obvious and concerning the symmetry we compute for $x \in (0,1)$
   \begin{align}
       u_t(1-x) & = u(\phi_t(1-x)) = u(1-x + t \tilde{\eta}(1-x)) 
       = u(1-x- t \tilde{\eta}(x)) \\ & = u(1- (x+t \tilde\eta(x))) = u(x+ t \tilde\eta(x)) = u(\phi_t(x))=u_t(x). 
   \end{align}
This and the minimalty of $u= u_0$ implies that $\frac{d}{dt} \mathcal{E}_\theta(u_t) \big\vert_{t = 0} = 0$. As a result, \eqref{eq:ddtequt} yields the claim.
\end{proof}

\subsection{Free boundary description}
\label{sec.Freeboundary}
The goal of this chapter is to show that the free boundary $\{u = 0\}$ is either empty or consists of precisely two points (which must due to symmetry be of the form $a, 1-a$ for some $a \in (0,\frac{1}{2})$). This will allow us to identify the Radon measure in  \eqref{eq:LemmaRadonMeasExplicit} in the latter case as a weighted sum of two Dirac measures, one located in $a$ and one located in $1-a$. Due to symmetry, the weights at both points must coincide. Hence, the Euler Lagrange equation reduces to a family of equations with two free parameters: One is the location of the Dirac measures $a$ and the other one is the weight. We also discuss the influence of both parameters on the shape of minimizers. It turns out that the weight is closely related to the gradient of minimizers at $a$, which motivates once again the study of the flatness question. 

 A first important observation that will be useful to describe the free boundary is that each minimizer is convex, also based on an idea in \cite[Lemma 2.1]{DDObstacle}.   
\begin{prop}\label{prop:convexity}
    Let $u$ be a minimizer of $\mathcal{E}_\theta$ in $\mathcal{A}_{\rm sym} (h)$. Then $u$ is convex. 
\end{prop}
\begin{proof}
    Let $u \in \mathcal{A}_{\rm sym}(h)$ be a minimizer. Now choose $w \in W^{1,2}(0,1)$ such that $w(0)=w(1)=h$ and 
    \begin{equation}
        \int_0^1 w'(x)^2 \; \mathrm{d}x = \inf \left\lbrace \int_0^1 f'(x)^2 \; \mathrm{d}x : f(0)=f(1)=h, f \leq u \right\rbrace.
    \end{equation}
    The existence of such $w$ follows from the standard direct methods in the calculus of variations. By strict convexity of the Dirichlet energy (and convexity of the admissible set), $w$ is also uniquely determined. Since $w$ and $w(1-\cdot)$ are both admissible and have the same Dirichlet energy, uniqueness of the minimizer yields $w = w(1-\cdot)$. Following the lines of \cite[Lemma 2.1]{DDObstacle} we can deduce that $w \in W^{2,2}(0,1)$ (in particular $w \in \mathcal{A}_{\rm sym}(h)$). Furthermore, (again following \cite[Lemma 2.1]{DDObstacle}) we have $w'' \geq 0$ a.e. and
    \begin{equation}\label{eq:w''}
        w'' = \begin{cases}
             0 & \textrm{a.e. on $\{ w < u \}$} \\ u'' & \textrm{a.e. on $\{ w = u \}$} 
        \end{cases}.
    \end{equation}
    Now we compute $\mathcal{E}_\theta(w)$. Since $w \leq u$ we have $|\{ w > 0 \}| \leq |\{ u > 0 \}|$. Moreover, \eqref{eq:w''} and the fact that $u'= w'$ a.e. on $\{u = w \}$ yield that 
    \begin{align}
        \mathcal{B}(w) = \int_0^1 \frac{w''^2}{(1+w'^2)^\frac{5}{4}} \; \mathrm{d}x = \int_{\{ w = u \}} \frac{w''^2}{(1+w'^2)^\frac{5}{4}} \; \mathrm{d}x = \int_{\{ w = u \}} \frac{u''^2}{(1+u'^2)^\frac{5}{4}} \; \mathrm{d}x  \leq \mathcal{B}(u). \label{eq:AbschwvsU}
    \end{align}
    Altogether we obtain $\mathcal{E}_\theta(w)\leq \mathcal{E}_\theta(u)$ and the inequality is strict unless $u'' = 0$ a.e. on $\{ u \neq w \}$ (see last estimate \eqref{eq:AbschwvsU}). A strict inequality contradicts minimality which is why we infer $u'' = 0$ a.e. on $\{ u  \neq w \}$. This implies that 
    \begin{equation}
        u'' = \begin{cases}
            w'' & \textrm{a.e. on $\{u = w\}$} \\ 0 & \textrm{a.e. on $\{u \neq w \}$} 
        \end{cases}.
    \end{equation}
    This and the fact that $w'' \geq 0$ a.e. yields $u'' \geq 0$ a.e. and hence convexity of $u$. 
\end{proof}

Notice that convex functions must either be monotone or decompose into a decreasing and an increasing function. Pure monotonicity is excluded due to the boundary values $u(0)=u(1) = h$ (unless $u$ is constant), hence the latter case applies. 

The interplay between this convexity and the scaling behavior of $\mathcal{B}$ paves the way for the announced non-flatness result for free boundary points, which we will prove in the subsequent proposition.

\begin{prop}\label{prop:Gradientnonzero}
     Let $u$ be a minimizer of $\mathcal{E}_\theta$ in $\mathcal{A}_{\rm sym} (h)$. Then $u' \neq 0$ on $\{u = 0 \}$.  
\end{prop}
\begin{proof}
Assume the opposite, i.e. $\{u =  u' = 0 \} \neq \emptyset$.
Define $x_1 := \min \{ x \in (0,1) : u(x) =u'(x) = 0 \}$ (attained since $\{u = u' = 0\}$ is a compact subset of $(0,1)$). Due to symmetry of $u$ we have that $x_1 \in (0, \frac{1}{2}]$. $C^1$-regularity and convexity of $u$ (see Proposition \ref{prop:convexity}) implies that for all $x \in (0,1)$ we have 
\begin{equation}\label{eq:upositiv}
    u(x) \geq u(x_1) + u'(x_1) (x-x_1) = 0 + 0 = 0.
\end{equation}
Hence $u \geq 0$. For this reason, everytime that $u$ attains the value $0$, $u$ attains its global minimum and thus $u' \vert_{\{u = 0 \} } = 0$, i.e. $\{ u = 0 \}= \{ u= u' = 0\}$.
\begin{flushleft}
    \textbf{Claim 1.} Under this assumption we have $\{ u = 0 \} = \{ u \leq 0 \} = [x_1,1-x_1].$
\end{flushleft}
That $\{ u = 0 \}=  \{ u \leq 0 \}$ follows from \eqref{eq:upositiv}. Now let $z \in [x_1,1-x_1]$, say $z = t x_1 + (1-t)(1-x_1)$ for some $t \in [0,1]$. As a consequence of \eqref{eq:upositiv} we have $u(z) \geq 0$. Moreover, we infer from convexity that 
\begin{equation}
    u(z) = u(tx_1 +(1-t)(1-x_1)) \leq t u(x_1) + (1-t)u(1-x_1)  = [t+ (1-t)]u(x_1) = 0. 
\end{equation}
This shows $[x_1,1-x_1] \subset \{u = 0\}$. Now let $z \in \{ u = 0 \} = \{ u= u' = 0 \}$. Then by definition of $x_1$ as a minimum we have $z \geq x_1$. Moreover $u(1-z) = 0$, i.e. $1-z \in \{ u = 0 \} = \{ u = u' = 0\}$, yielding also $1-z \geq x_1$, that is $z \leq 1-x_1$. We have now also shown that $\{ u = 0 \} \subset [x_1,1-x_1]$ and infer the claim. An immediate consequence is that 
\begin{equation}\label{eq:energycomp}
    \mathcal{E}_\theta(u) = \mathcal{B}(u) + 2x_1 \theta.
\end{equation}
\begin{flushleft}
    \textbf{Claim 2.} $x_1< \frac{1}{2}$ 
\end{flushleft}
We assume the opposite, i.e. $x_1 = \frac{1}{2}$. We observe from \eqref{eq:energycomp} that
$
    \mathcal{E}_\theta(u) \geq \mathcal{B}(u) + \theta. 
$
Since however by Lemma \ref{prop:theta} we have $\mathcal{E}_\theta(u) \leq \theta$ we infer $\mathcal{B}(u)=0$. This is only possible if $u'' =0 $ which yields (together with $u(0) = u(1) = h$) that $u \equiv h$. But then $\{u = 0\} = \emptyset$ which contradicts our assumption that $\{u = u' = 0 \} \neq \emptyset.$  Claim 2 is shown. 

Now choose some $\lambda > 1$ such that $\lambda x_1< \frac{1}{2}$ and define a comparison function $u_\lambda: [0,1] \to \mathbb{R}$ via the following scaling operation (see also Figure \ref{fig:placeholder3})
\begin{equation}
    u_\lambda(x):= \begin{cases}
         h + \lambda (u(\frac{x}{\lambda}) - h) & x \in [0,\lambda x_1] \\ (1-\lambda) h & x \in [\lambda x_1, \frac{1}{2}] \\
         \textrm{extended in a symmetric way} & x \in [\frac{1}{2},1]
    \end{cases}. 
\end{equation}
\begin{figure}[h!]
    \centering
    \includegraphics[width=0.25\linewidth]{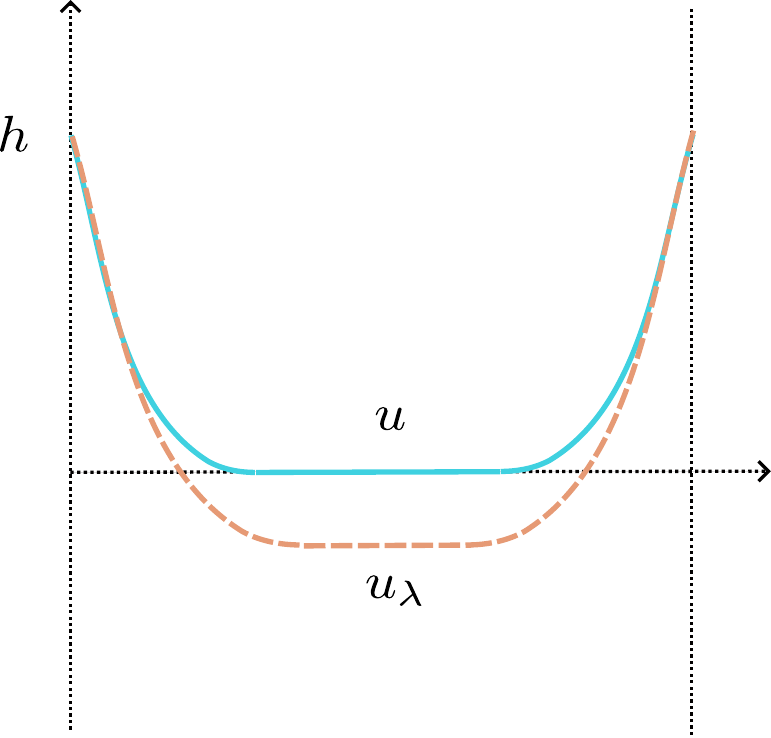}
    \caption{The comparison function $u_\lambda$ for some $\lambda >1$.}
    \label{fig:placeholder3}
\end{figure}

Since $u_\lambda(0) = h+ \lambda(h-h)= h$  we infer that $u_\lambda \in \mathcal{A}_{\rm sym}(h)$. Due to the fact that for all $x \in [0,\lambda x_1]$ we have $u_\lambda'(x) = u'(\frac{x}{\lambda})$ and $u_\lambda''(x) =  \frac{1}{\lambda}u''(\frac{x}{\lambda})$ we obtain
\begin{align}
    \mathcal{B}(u_\lambda)  & = \frac{2}{\lambda^2} \int_0^{\lambda x_1} \frac{u''(\frac{x}{\lambda})^2}{(1+u'(\frac{x}{\lambda})^2)^\frac{5}{2}} \; \mathrm{d}x = \frac{2}{\lambda} \int_0^{x_1} \frac{u''(x)^2}{(1+u'(x)^2)^\frac{5}{2}} \; \mathrm{d}x  \\ & = \frac{1}{\lambda} \int_{(0,x_1) \cup(1-x_1,1)}  \frac{u''(x)^2}{(1+u'(x)^2)^\frac{5}{2}} \; \mathrm{d}x \leq \frac{1}{\lambda} \mathcal{B}(u). \label{eq:togetherwith}
\end{align}
\begin{flushleft}
    \textbf{Claim 3.} There exists $\lambda_0 > 1$ such that for all $\lambda \in (1,\lambda_0)$ we have $\{u_\lambda \leq 0 \} \supset [x_1,1-x_1]$.
\end{flushleft}
Clearly, for each $\lambda > 1$ we have $u_\lambda\vert_{[\lambda x_1,\frac{1}{2}]} \equiv h (1-\lambda) <0$. Therefore (and by symmetry) it suffices to show that $u_\lambda\vert_{[x_1,\lambda x_1]} \leq 0$ for $\lambda$ sufficiently close to $1$. Now suppose that there exists some $x \in [x_1,\lambda x_1]$ such that $0 \leq u_\lambda(x) = h + \lambda (u(\frac{x}{\lambda})- h)$. This implies 
\begin{align}
    h (\lambda -1)  & \leq \lambda u(\tfrac{x}{\lambda}) = \lambda [u(\tfrac{x}{\lambda}) - u(x_1)]  \leq  \lambda \sup_{z \in [\frac{x_1}{\lambda},x_1]} |u'(z)|  \cdot |\tfrac{x}{\lambda} - x_1| \\ &  \leq   \sup_{z \in [\frac{x_1}{\lambda},x_1]} |u'(z)|  \cdot |x- \lambda x_1| \leq  \sup_{z \in [\frac{x_1}{\lambda},x_1]} |u'(z)|   x_1 (\lambda-1) \leq  \tfrac{1}{2}\sup_{z \in [\frac{x_1}{\lambda},x_1]} |u'(z)|  (\lambda-1). 
\end{align}
Dividing by $\lambda - 1$ we find $h \leq \tfrac{1}{2}\sup_{z \in [\frac{x_1}{\lambda},x_1]} |u'(z)|$. However, since $u'(x_1)= 0$ and $u'$ is continuous at $x_1$ we can find $\lambda_0 > 1$ sufficiently close to $1$ such that for all $\lambda \in (1, \lambda_0)$ one has the opposite inequality $h > \tfrac{1}{2}\sup_{z \in [\frac{x_1}{\lambda},x_1]} |u'(z)|$. The contrapositive of the above implication shows that then there can not exist $x \in [x_1,\lambda x_1]$ such that $u_\lambda(x) \geq 0$. In particular, $[x_1,\lambda x_1] \subset \{ u_\lambda \leq 0 \}$ and Claim 3 follows. 



From Claim 3 we infer that for all $\lambda \in (1,\lambda_0)$ we have 
\begin{equation}
    |\{ u_\lambda  > 0 \}| \leq 2x_1 = |\{u > 0 \}|.
\end{equation}
This and \eqref{eq:togetherwith} imply that for $\lambda \in (1,\lambda_0)$ we have
\begin{equation}
    \mathcal{E}_\theta(u_\lambda) \leq  \frac{1}{\lambda} \mathcal{B}(u) + \theta |\{u > 0 \}|= \mathcal{E}_\theta(u) + \left( \frac{1}{\lambda}- 1 \right) \mathcal{B}(u).  
\end{equation}
Since $u$ is a minimizer we must also have $\mathcal{E}_\theta(u_\lambda) \geq \mathcal{E}_\theta(u)$ and thus $(\frac{1}{\lambda} - 1) \mathcal{B}(u) \geq 0$. This can only work if $\mathcal{B}(u)= 0$, i.e. $u'' \equiv 0$ a.e. But then the fact that $u \in \mathcal{A}_{\rm sym}(h)$ yields $u \equiv h$, contradicting once again the assumption that $\{ u = u' = 0\} \neq \emptyset$.
\end{proof}


\begin{prop}\label{eq:prop311}
    Let $u$ be a minimizer of $\mathcal{E}_\theta$ in $\mathcal{A}_{\rm sym} (h)$. Then either $\{u = 0 \} = \emptyset$ or
    there exists  $a \in (0,\frac{1}{2})$ such that $\{u = 0\} = \{ a, 1-a\}$. In the latter case we have $u> 0$ on $(0,a) \cup (1-a,1)$ and $u<0$ on $(a,1-a)$.
\end{prop}
\begin{proof}
Suppose that $\{u = 0 \} \neq \emptyset$.
    Due to the previous proposition we have $\{ u= 0 \} = \{ u =0, u' \neq 0\}$. The set on the right hand side is obviously discrete and the set on the left hand side is compact. Since discrete and compact sets are finite we must have $\{u = 0 \} = \{ a_1,...,a_m\}$ for some $a_1,...,a_m \in (0,1)$ and $m \in \mathbb{N}$. Notice that  symmetry implies $u'(\frac{1}{2}) =0 $ and hence $\frac{1}{2} \not \in \{ u = 0, u' \neq 0 \}= \{ u = 0\}$. Moreover, if $b \in \{u = 0\}$ then also $1-b \in \{ u = 0 \}$. Both previous observations imply that $m \in \mathbb{N}$ is even. It remains to show that $m = 2$. Suppose that $m \neq 2$. Evenness yields $m \geq 4$. Let $a_1< a_2 < a_3 < a_4$ be the smallest four elements of $\{u=0\}$. In particular $u(a_i) = 0$ and $u'(a_i) \neq 0$ for all $i \in \{1,2,3,4\}$. As a consequence of the nonvanishing gradient and the intermediate value theorem we have $u\vert_{(0,a_1)} > 0, u\vert_{(a_1,a_2)}< 0$ and $u\vert_{(a_2,a_3)} > 0$. 
    In particular, $u(a_2)=u(a_3) =0$ and $u ( \frac{a_2+a_3}{2}) > 0$. Convexity of $u$ however yields 
    \begin{equation}
        u ( \tfrac{a_2+a_3}{2}) \leq \tfrac{1}{2}(u(a_2)+ u(a_3))= \tfrac{1}{2}(0+0) = 0, 
    \end{equation}
    a contradiction. We infer that $m= 2$, i.e. $\{u = 0 \}= \{ a_1,a_2\}$ for some $a_1,a_2 \in (0,1)$ such that $a_1< a_2$. Symmetry of $u$ yields that $a_1 \in (0,\frac{1}{2})$ and $a_2 = 1-a_1$ and hence implies that $\{u = 0\} = \{ a,1-a \}$ with $a:= a_1$. The last sentence of the statement follows directly from Proposition \ref{prop:Gradientnonzero}, as the nonvanishing gradient enforces a sign change at $a$ and $1-a$.  
\end{proof}

\begin{remark}\label{rem:0indermitte}
    We know from Proposition \ref{eq:prop311} that the free boundary consists of two points $\{a,1-a\}$ for some $a \in (0,\frac{1}{2})$, which divides the set $\{ u \neq 0 \}$ into three components, namely $(0,a),(a,1-a)$ and $(1-a,1)$ on all of which $u$ has a strict sign. Proposition \ref{prop:preevious} yields that $u$ is smooth on all of there components and $\frac{d}{dx} \frac{v_u'}{(1+u'^2)^\frac{5}{4}} = 0$. This implies that we have $C_1,C_2,C_3 \in \mathbb{R}$ such that 
    \begin{equation}\label{eq:evalaux}
        \frac{v_u'(x)}{(1+u'(x)^2)^\frac{5}{4}} = \begin{cases}
            C_1 & x \in (0,a) \\ C_2 & x \in (a,1-a) \\ C_3 & x \in (1-a,1)
        \end{cases}.
    \end{equation}
Note that $u= u(1-\cdot)$ yields $u'^2= u'(1-\cdot)^2$,  $v_u = v_u(1-\cdot)$ and $v_u' = -v_u'(1-\cdot)$. In particular, $\frac{v_u'}{(1+u'^2)^\frac{5}{4}}  = - \frac{v_u'}{(1+u'^2)^\frac{5}{4}} (1- \cdot)$. This antisymmetry and \eqref{eq:evalaux} immediately yields $C_3 = -C_1$. We claim also that $C_2 = 0$. Indeed, consider some $b \in (a,\frac{1}{2})$. Then $1-b \in (\frac{1}{2},1-a)$ and hence 
\begin{equation}
     \frac{v_u'(b)}{(1+u'(b)^2)^\frac{5}{4}} = C_2 = \frac{v_u'(1-b)}{(1+u'(1-b)^2)^\frac{5}{4}}.
\end{equation}
Using again the antisymmetry we obtain  $\frac{v_u'(1-b)}{(1+u'(1-b)^2)^\frac{5}{4}} = -  \frac{v_u'(b)}{(1+u'(b)^2)^\frac{5}{4}}$ which yields $C_2= 0$. As a result we have for some $C_1 \in \mathbb{R}$
\begin{equation}\label{eq:LHSThirsORDER}
    \frac{v_u'(x)}{(1+u'(x)^2)^\frac{5}{4}} = \begin{cases}
        C_1 & x \in (0,a) \\ 0 & x \in (a,1-a) \\ -C_1 & x \in (1-a,1)
    \end{cases}.
\end{equation}
We will use this to characterize $u$ more explicitly.
\end{remark}
The zero in the middle part of \eqref{eq:LHSThirsORDER} allows us to give an explicit shape to minimizers on $(a,1-a)$.

\begin{prop}\label{prop:cshapeparameter}
     Let $u$ be a minimizer of $\mathcal{E}_\theta$ in $\mathcal{A}_{\rm sym} (h)$ with $\{u = 0 \} = \{ a, 1-a\}$ for some $a \in (0,\frac{1}{2})$. Then there exists $c \in (0,c_0)$ such that for all $x \in (a,1-a)$ we have 
     \begin{equation}\label{eq:uina1-a}
         u(x) = -(1-2a) u_{c}( \tfrac{x-a}{1-2a}), 
     \end{equation}
     where $u_c$ is defined as in \eqref{eq:uCAUx}.
     In particular, we have $u'(a)  = -G^{-1}(\frac{c}{2})$ and $v_u(a) = \frac{c}{1-2a}$ (i.e. $u''(a)= (1+ G^{-1}(\frac{c}{2})^2)^\frac{5}{4} \frac{c}{1-2a}$).
\end{prop}
\begin{proof}
    As argued in Remark \ref{rem:0indermitte} we have $\frac{v_u'}{(1+u'^2)^\frac{5}{4}}= 0$ on $(a,1-a)$. This implies that $v_u$ is constant on $(a,1-a)$. Let $\tilde{c} \in \mathbb{R}$ be such that $v_u(x) = \tilde{c}$ for all $x \in (a,1-a)$. Notice that $\tilde{c}\geq 0$ due to convexity of $u$, cf. Proposition \ref{prop:convexity}.   
    Since $v_u= (G \circ u')'$ there exists $d \in \mathbb{R}$ such that $G(u'(x)) ) = \tilde{c}x+ d$. The fact that $u'(\frac{1}{2}) = 0$ yields $d= - \frac{\tilde{c}}{2}$ and $u'(x) = G^{-1}(\tilde{c}x - \frac{\tilde{c}}{2})  = - G^{-1}( \frac{\tilde{c}}{2} - \tilde{c}x)$ for all $x \in (a,1-a)$.
    Note that $-\infty < u'(a) =- G^{-1}(\frac{1}{2}\tilde{c}(2a-1))$, which implies that $\tilde{c}< \frac{c_0}{1-2a}$. Now define $c:= (1-2a)\tilde{c} \in [0,c_0)$ and $w: [a,1-a] \to \mathbb{R}$ via $w(x) := -(1-2a) u_c(\frac{x-a}{1-2a})$. Observe that for $x \in (a,1-a)$ we have
    \begin{equation}
        w'(x) = -u_c'(\tfrac{x-a}{1-2a}) = -G^{-1}(\tfrac{c}{2}- c \tfrac{x-a}{1-2a}) = -G^{-1} ( \tfrac{1}{2} \tfrac{c}{1-2a} - \tfrac{c}{1-2a}x) = -G^{-1}(\tfrac{\tilde{c}}{2}-\tilde{c}x) = u'(x).  
    \end{equation}
    This and the fact that $w(a)= 0= u(a)$ imply that $u \equiv w$ on $(a,1-a)$. We infer \eqref{eq:uina1-a}. Due to the fact that $u,u',u''$ are continuous at $a$ (see Lemma \ref{lem:RadonMeas}) we may use \eqref{eq:uina1-a} to compute the asserted formulas  $u'(a)= - u_c'(0) = G^{-1}(\frac{c}{2})$ and $v_u(a) = \frac{c}{1-2a}$. Notice that the formula for $u'(a)$ and Proposition \ref{prop:Gradientnonzero} imply that $c \neq 0$ and hence $c \in (0,c_0)$  
\end{proof}

\begin{remark}\label{rem:v_ugleichconstant}
    Proposition \ref{prop:cshapeparameter} and \eqref{eq:LHSThirsORDER} actually imply that $v_u \equiv \frac{c}{1-2a}$ on $[a,1-a].$
\end{remark}
With our previous insights we can characterize the Radon measure $\mu$ from Lemma \ref{lem:RadonMeas} more explicitly and express the constant $C_1$  in \eqref{eq:LHSThirsORDER} in terms of the new shape parameter $c$.
\begin{corollary}\label{cor:314}
    Suppose that $u \in \mathcal{A}_{\rm sym}(h)$ is a minimizer of $\mathcal{E}_\theta$. Let $a \in (0,\frac{1}{2})$ be such that $\{u = 0\} = \{ a,1-a \}$ and $c \in (0,c_0)$ be as in Proposition \ref{prop:cshapeparameter}. Then we have $c > 0$ and 
    \begin{equation}\label{vuprimme}
        \frac{v_u'(x)}{(1+u'(x)^2)^\frac{5}{4}} = \begin{cases}
            \frac{\theta}{2G^{-1}(\frac{c}{2})} & x \in (0,a) \\
            0 & x \in (a,1-a) \\
            -\frac{\theta}{2G^{-1}(\frac{c}{2})} & x \in (1-a,1)
        \end{cases}.
    \end{equation}
\end{corollary}
\begin{proof}
    Let $\mu$ be the Radon measure in Lemma \ref{lem:RadonMeas}. Note that $\mu$ is supported on $\{a,1-a\}$ and hence must be of the form $\lambda_1 \delta_a + \lambda_2 \delta_{1-a}$, where $\lambda_1,\lambda_2 \in \mathbb{R}$ and $\delta_x$ denotes the Dirac measure centered at $x$.  The coefficients $\lambda_1,\lambda_2$ can be determined comparing both sides of Lemma \ref{eq:LemmaRadonMeasExplicit}. Indeed, for any $\eta \in C_0^\infty((0,1))$ we now have with Proposition \ref{prop:cshapeparameter}
    \begin{equation}
        \int_{\{ u = 0 \} } u' \eta \; \mathrm{d}\mu = \lambda_1 u'(a) \eta(a) + \lambda_2u'(1-a) \eta(1-a) =  -\lambda_1 G^{-1}(\tfrac{c}{2})  \eta(a) + \lambda_2 G^{-1}(\tfrac{c}{2}) \eta(1-a)
    \end{equation}
    and
    \begin{equation}
        - \theta \int_{\{u > 0 \}} \eta' \; \mathrm{d}x= -\theta \int_0^a \eta' \; \mathrm{d}x  - \theta \int_{1-a}^1 \eta' \; \mathrm{d}x = - \theta\eta(a) + \theta \eta(1-a).
    \end{equation}
    This implies that $\lambda_1 = \lambda_2 =  \frac{\theta}{G^{-1}(\frac{c}{2})}$, i.e. $\mu = \frac{\theta}{G^{-1}(\frac{c}{2})}(\delta_a + \delta_{1-a})$. Having found this representation, we can compute both sides of \eqref{eq:todo} and infer the value of $C_1$ in \eqref{eq:LHSThirsORDER}. Let $\varphi \in C_0^\infty((0,1))$. For the left hand side of \eqref{eq:todo} we compute with Lemma \ref{lem:33} and \eqref{eq:LHSThirsORDER}
    \begin{equation}
        D \mathcal{B}(u) (\varphi) =  - 2 \int_0^1 \frac{v_u'(x)}{(1+ u'(x)^2)^\frac{5}{4}} \varphi'(x) = -2 \int C_1(\chi_{(0,a)}(x)- \chi_{(1-a,1)}(x)) \varphi'(x) \; \mathrm{d}x = - 2C_1 (\varphi(a) +  \varphi(1-a)) 
    \end{equation}
    For the right hand side of \eqref{eq:todo} we have by the formula for $\mu$
    \begin{equation}
        - \int \varphi \; \mathrm{d}\mu = - \frac{\theta}{G^{-1}(\frac{c}{2})}(\varphi(a) + \varphi(1-a) ).
    \end{equation}
    Since $\varphi$ was arbitrary we obtain $C_1 = \frac{\theta}{2G^{-1}(\frac{c}{2})}$ and the claim follows with \eqref{eq:LHSThirsORDER}.
\end{proof}

\begin{remark}\label{rem:nitchW41}
Note that Corollary \ref{cor:314} implies that the third derivative $u'''$ must have two jumps at $a$ and $1-a$. This in particular shows that $u \not \in W^{4,1}(0,1)$ and hence implies optimality of the $W^{3,\infty}$-regularity as asserted in Lemma \ref{lem:RadonMeas}.    
\end{remark}

A corollary which will be of importance later is the \textit{strict convexity} of $u$ on $(0,1)$, further improving Proposition \ref{prop:convexity}. 

\begin{corollary}[Strict Convexity]\label{cor:strictconvexity}
    Suppose that $u$ is a minimizer of $\mathcal{E}_\theta$ in $\mathcal{A}_{\rm sym}(h)$ with $\{u = 0 \} \neq \emptyset$. Then $u'' > 0$ on $(0,1)$.
\end{corollary}
\begin{proof}
    By symmetry, it suffices to show the claim on $(0,\frac{1}{2}]$. Therefore let $x \in (0,\frac{1}{2}]$. We distinguish two cases, namely $x \in (0,a)$ and $x \in [a, \frac{1}{2}]$. First let $x \in (0,a)$. We deduce from \eqref{vuprimme} that 
    $v_u'  = (1+u'^2)^\frac{5}{4} G^{-1}(\frac{c}{2})^{-1}$ lies in $C^0((0,a))$ and moreover 
    Lemma \ref{lem:RadonMeas} implies that $v_u \in C^0([0,a])$. As a result we have $v_u \in C^1((0,a)) \cap C^0([0,a])$. This makes the mean value theorem applicable and we deduce that 
    \begin{equation}
        v_u(x) =  v_u(0) +  v_u'(\xi)(x-0)   \qquad \textrm{for some $\xi \in (0,x)$.}
    \end{equation}
    Using that by Lemma \ref{lem:36} we have $v_u(0) = 0$ and by \eqref{vuprimme} we have $v_u'(\xi) =(1+u'(\xi)^2)^\frac{5}{4} G^{-1}(\frac{c}{2})^{-1} > 0$ we obtain $v_u(x) > 0$ and as a consequence $u''(x) > 0$. Next let $x \in [a, \frac{1}{2}]$. Then \eqref{eq:uina1-a} and \eqref{v_uisconstant} yield 
    \begin{equation}
        u''(x) = -\frac{1}{(1-2a)} u_c''(\tfrac{x-a}{1-2a}) = \frac{1}{(1-2a)} c (1+ u_c'(\tfrac{x-a}{1-2a})^2)^\frac{5}{4} > 0. 
    \end{equation}
    The claim follows.
\end{proof}

\begin{proof}[Proof of Theorem \ref{thm3}.]
Let $u \in \mathcal{A}_{\rm sym}(h)$ be a minimizer and $u \not \equiv h$. If $u$ does not change sign then
\begin{equation}
    \mathcal{E}_\theta(u) = \mathcal{B}(u) + \theta |(0,1)| = \mathcal{B}(u) + \theta \geq \theta = \mathcal{E}_\theta(h) \geq \mathcal{E}_\theta(u), 
\end{equation}
where the last estimate is due to the minimality of $u$. Hence $\mathcal{B}(u)= 0$, resulting in  the contradiction $u \equiv h$. As a consequence $u$ is a sign-changing minimizer, i.e. $\{u = 0 \} \neq \emptyset.$
    The following list summarizes where the proofs of all the needed assertions sign-changing minimizers $u \in \mathcal{A}_{\rm sym}(h)$ of $\mathcal{E}_\theta$ can be found 
    \begin{itemize}[itemsep=-0.5pt]
        \item $W^{3,\infty}$-regularity: Lemma \ref{lem:RadonMeas}.
        \item Boundary conditions $u''(0)= u''(1)= 0$: Lemma \ref{lem:36}.
        \item Strict convexity on $(0,1)$: Corollary \ref{cor:strictconvexity}.
        \item Existence of $a \in (0,1)$ such that $\{u = 0\}= \{ a, 1-a\}$: Proposition \ref{eq:prop311}.
        \item Nonvanishing gradient at the free boundary, i.e. $u'(a),u'(1-a) \neq 0$: Proposition \ref{prop:Gradientnonzero}.
        \item Three-elastica-coupling: Proposition \ref{prop:preevious} (using that by Proposition \ref{eq:prop311} $u$ does not have a zero on $(0,a)$, on $(a,1-a)$ and on $(1-a,1)$). Strict convexity yields that none of them are lines. \qedhere
    \end{itemize}
\end{proof}
We have now reduced the Euler Lagrange equation to a two-parameter family with parameters $a$ and $c$, cf. Corollary \ref{cor:314}. Due to the nonlinearity and higher order of the Euler-Lagrange equation is is not clear whether the two parameters determine a unique candidate for a minimizer. We shed some light on this in the following remark. 
\begin{remark}\label{ref:rem3.14}
Suppose that $u \in \mathcal{A}_{\rm sym}(h)$ is a minimizer of $\mathcal{E}_\theta$ and $\{u = 0 \}\neq \emptyset$. Let $a \in (0,\frac{1}{2})$ be such that $\{u = 0\} = \{ a,1-a \}$ and $c \in [0,c_0)$ be as in Proposition \ref{prop:cshapeparameter}. 
By Proposition \ref{prop:cshapeparameter} we have the boundary conditions  $u(a) = 0, u'(a) = -G^{-1}(\frac{c}{2}),u''(a)=  (1+ G^{-1}(\frac{c}{2})^2)^\frac{5}{4} \frac{c}{1-2a}$. Using this and \eqref{vuprimme} $u\vert_{[0,a]}$ is determined uniquely by the third order initial value problem
\begin{equation}\label{eq:uniquesol}
       \begin{cases} \frac{d}{dx} \frac{u''(x)}{(1+u'(x)^2)^\frac{5}{4}} = \frac{\theta}{2G^{-1}(\frac{c}{2})} (1+ u'(x)^2)^\frac{5}{4}  & x \in [0,a] \\ 
       u(a) = 0 \\ u'(a) = -G^{-1}(\frac{c}{2})  \\ u''(a) =  (1+ G^{-1}(\frac{c}{2})^2)^\frac{5}{4} \frac{c}{1-2a}
       \end{cases}.
    \end{equation}
Moreover, also $u\vert_{[a,1-a]}$ is determined uniquely by $a$ and $c$ due to Proposition \ref{prop:cshapeparameter}. Furthermore, $u\vert_{[1-a,1]}$ is determined uniquely by $u\vert_{[0,a]}$ due to symmetry. This shows that once $a$ and $c$ are known the minimizer is determined uniquely. In particular, each pair of values of $a$ and $c$ gives a unique candidate for a minimizer. There is actually a subtle point about this statement: We only obtain a candidate if the solution of \eqref{eq:uniquesol} exists on the whole interval $[0,a]$. This is nontrivial for a nonlinear ODE, a priori we only obtain existence on a maximal existence interval $(x_-,x_+)$ containing $a$ and it is not immediately clear whether $x_-< 0$. Taking into account this detail we can only say that for each $a \in (0,\frac{1}{2}),c \in (0,c_0)$ we have \textit{at most} one candidate for the solution. This will be sufficient for the arguments to come. We will address the role of this subtility when needed. 
\end{remark}

We conclude this section by a useful formula for the boundary derivative $u'(0)$ of minimizers $u \in \mathcal{A}_{\rm sym}(h)$ that will be of use later on.

\begin{prop}\label{eq:energyformula}
    Let $u \in \mathcal{A}_{\rm sym}(h)$ be a minimizer of $\mathcal{E}_\theta$. Let $a \in (0,\frac{1}{2})$ be such that $\{u = 0\} = \{ a,1-a \}$ and $c \in (0,c_0)$ be as in Proposition \ref{prop:cshapeparameter}.  Then  $u'(0) = - G^{-1}(\tfrac{c}{2}) (1+ \frac{1}{\theta}\tfrac{c^2}{(1-2a)^2}).$
\end{prop}
\begin{proof}
    
    Recall from Corollary \ref{cor:314} that for all $x \in (0,a)$ one has 
    \begin{equation}
        v_u'(x) = \frac{\theta}{2 G^{-1}(\frac{c}{2})} (1+u'(x)^2)^\frac{5}{4}.
    \end{equation}
    We multiply by $2v_u(x)$ and use the definition of $v_u$ on the right hand side to find 
   \begin{equation}
        2 v_u(x) v_u'(x) = \frac{\theta}{G^{-1}(\frac{c}{2})} u''(x).
    \end{equation}
    Integrating from $a$ to $x$ we find 
   $
        v_u(x)^2 - v_u(a)^2 = \frac{\theta}{G^{-1}(\frac{c}{2})} (u'(x) - u'(a)).
  $
    Using that by Proposition \ref{prop:cshapeparameter} we have $v_u(a)= \frac{c}{1-2a}$ and $u'(a)= - G^{-1}(\frac{c}{2})$ we infer that for all $x \in (0,a)$
    \begin{equation}\label{eq:Bendingdarstellung}
        v_u(x)^2= \frac{c^2}{(1-2a)^2} +   \frac{\theta}{G^{-1}(\frac{c}{2})} (u'(x) + G^{-1}(\tfrac{c}{2})) = \frac{c^2}{(1-2a)^2} + \theta + \frac{\theta u'(x)}{G^{-1}(\frac{c}{2})}.
    \end{equation}
   Plugging in $x= 0 $  and using that by Lemma \ref{lem:36} $v_u(0)=0$ we find the formula for $u'(0)$ given in the statement.
\end{proof}


\section{The improved nonexistence threshold}\label{sec:4}
We have now reduced the candidates for sign-changing minimizers to a two-parameter family with two parameters $a \in (0, \frac{1}{2})$ and $c \in (0,c_0)$. These two parameters give rise to at most one candidate $u_{a,c}$ for a minimizer  by using the formula from  Proposition \ref{prop:cshapeparameter} on $[a,1-a]$ and (if possible) solving \eqref{eq:uniquesol} on $[0,a]$ (and then extending symmetrically). However, the parameters $a$ and $c$ are not independent as this procedure ignores the boundary condition $u(0)= h$. We will next discuss how one has to choose $a$ and $c$ to obtain a prescribed height $h$. We will also see that the map $(a,c) \mapsto u_{a,c}(0)$ is not surjective and hence not all values of $h$ can be reached. This will result in the improved nonexistence threshold in Theorem \ref{thm2}. For notational convenience we introduce new parameters $b$ and $s$ which can be computed from $a$ and $c$. These are defined as follows
\begin{equation}\label{eq:NewParameters}
    b := \tfrac{c^2}{(1-2a)^2}, \qquad s := G^{-1}(\tfrac{c}{2}) \left( 1+ \tfrac{1}{\theta} \tfrac{c^2}{(1-2a)^2} \right).
\end{equation}
Notice that these have an interpretation, namely according to Proposition  \ref{prop:cshapeparameter} we have that $b = v_{u_{a,c}}(a)^2$ and according to Proposition \ref{eq:energyformula} we have $s = |u_{a,c}'(0)|$. Also observe that $(a,c)$ can be retrieved uniquely from $(b,s)$ via $c= 2 G(\frac{\theta s}{\theta+b})$ and $a = \frac{1}{2}- \frac{1}{\sqrt{b}}G(\frac{\theta s}{\theta+b})$.

\begin{prop}\label{Prop:4.1}
    Let $u \in \mathcal{A}_{\rm sym}(h)$ be a minimizer with $\{u = 0 \} \neq \emptyset$.  Let $a \in (0,\frac{1}{2})$ be such that $\{u = 0\} = \{ a,1-a \}$ and $c \in (0,c_0)$ be as in Proposition \ref{prop:cshapeparameter}. Define $b$ and $s$ as in \eqref{eq:NewParameters}.
    Then 
    \begin{equation}\label{eq:monotoneinb}
        h =  \frac{1}{\sqrt{\theta+b}}\int_{\theta/(\theta+b)}^1 \frac{s^2y}{(1+s^2y^2)^\frac{5}{4}} \frac{1}{\sqrt{1-y}} \; \mathrm{d}y. 
    \end{equation}
\end{prop}
\begin{proof}
    Following the lines of the derivation of  \eqref{eq:Bendingdarstellung} in the proof of Proposition \ref{eq:energyformula} we derive that for all $x \in (0,a)$ we have 
    \begin{equation}
        v_u(x)^2 = \frac{c^2}{(1-2a)^2} + \theta +  \theta\frac{u'(x)}{G^{-1}(\frac{c}{2})} = \theta\frac{s + u'(x)}{G^{-1}(\frac{c}{2})}.
    \end{equation}
    Notice that due to convexity and the fact that $u'(\frac{1}{2}) = 0, u'(0)=-s$ (see Proposition \ref{eq:energyformula})  we have $-s \leq u'(x) \leq 0$ for all $x \in (0,a]$ and $u''(x) \geq 0$ for all $x \in (0,a]$. Using this and $v_u(x)= \frac{u''(x)}{(1+ u'(x)^2)^\frac{5}{4}}$ we obtain 
    \begin{equation}
        u''(x) = \sqrt{\theta}(1+u'(x)^2)^\frac{5}{4}\frac{\sqrt{s + u'(x)}}{G^{-1}(\frac{c}{2})^\frac{1}{2}} \qquad \textrm{for all $x \in (0,a].$}
    \end{equation}
    Next define $F: [-s, \infty) \to \mathbb{R}$
    \begin{equation}
        F(r)  := \int_{-s}^r \frac{1}{(1+t^2)^\frac{5}{4}} \frac{1}{\sqrt{s+t}} \; \mathrm{d}t. 
    \end{equation}
    Notice that $F$ is well-defined since the integrand is positive and
    \begin{equation}
        \int_{-s}^r \frac{1}{(1+t^2)^\frac{5}{4}} \frac{1}{\sqrt{s+t}} \; \mathrm{d}t \leq \int_{-s}^r \frac{1}{\sqrt{s+t}} \; \mathrm{d}t \leq \left[ 2 \sqrt{s+t} \right]_{-s}^r = 2 \sqrt{s+r}. 
    \end{equation}
    This computation also shows that $F$ is continuous at $r = -s$ and attains the value $F(-s)= 0$. Moreover one readily checks that $F$ is differentiable on $(-s,\infty)$, with $F'(r)< 0$ for all $r \in (-s,\infty)$. 
    Due to  Corollary \ref{cor:strictconvexity} we have $u'(x) \in (-s, \infty)$ for all $x \in (0,a)$. As a result we can apply the chain rule and find for each $x \in (0,a)$
    \begin{equation}
        \frac{d}{dx} (F \circ u')(x) = F'(u'(x)) u''(x) = \frac{u''(x)}{(1+u'(x)^2)^\frac{5}{4} } \frac{1}{\sqrt{s+u'(x)}} = \frac{\sqrt{\theta}}{G^{-1}(\frac{c}{2})^\frac{1}{2}}.
    \end{equation}
    Using this, the continuity of $F \circ u'$ at $0$ and the fact that that $(F \circ u')(0)= F(u'(0))=F(-s) = 0$ we find 
    \begin{equation}
        F(u'(x))= \frac{\sqrt{\theta}}{G^{-1}(\frac{c}{2})^\frac{1}{2}} x \qquad \textrm{for all $x \in [0,a]$. }
    \end{equation}
   Standard results on inverse functions yield that $F^{-1}$ is continuous on $[0, \frac{\sqrt{\theta}}{G^{-1}(\frac{c}{2})^{1/2}}a ]$ and differentiable in $(0,  \frac{\sqrt{\theta}}{G^{-1}(\frac{c}{2})^{1/2}}a ]$ and 
   \begin{equation}
       u'(x) = F^{-1} \left(  \frac{\sqrt{\theta}}{G^{-1}(\frac{c}{2})^\frac{1}{2}} x \right) \qquad \textrm{for all $x \in [0,a]$.}
   \end{equation}
We obtain 
\begin{align}
    h & = u(0) = u(0) - u(a) = \int_a^0 u'(t) \; \mathrm{d}t = \int_a^0 F^{-1} \left(  \frac{\sqrt{\theta}}{G^{-1}(\frac{c}{2})^\frac{1}{2}} t \right)  \; \mathrm{d}t= \frac{G^{-1}(\tfrac{c}{2})^\frac{1}{2}}{\sqrt{\theta}}\int_{a \sqrt{\theta} G^{-1}(\frac{c}{2})^{-\frac{1}{2}}}^0 F^{-1}(y) \, \mathrm{d}y
    \\ & = \frac{G^{-1}(\tfrac{c}{2})^\frac{1}{2}}{\sqrt{\theta}}\int_{F(u'(a))}^{F(u'(0))} F^{-1}(y) \, \mathrm{d}y = \frac{G^{-1}(\tfrac{c}{2})^\frac{1}{2}}{\sqrt{\theta}} \int_{u'(a)}^{u'(0)} F^{-1}(F(y)) F'(y) \; \mathrm{d}y
    \\ & =\frac{G^{-1}(\tfrac{c}{2})^\frac{1}{2}}{\sqrt{\theta}}\int_{u'(a)}^{u'(0)} y \frac{1}{(1+y^2)^\frac{5}{4}} \frac{1}{\sqrt{s+y}} \; \mathrm{d}y = \frac{G^{-1}(\tfrac{c}{2})^\frac{1}{2}}{\sqrt{\theta}}\int_{-G^{-1}(\frac{c}{2})}^{-s} y \frac{1}{(1+y^2)^\frac{5}{4}} \frac{1}{\sqrt{s+y}} \; \mathrm{d}y  \\ & = - \frac{G^{-1}(\tfrac{c}{2})^\frac{1}{2}}{\sqrt{\theta}}\int_{-s}^{-G^{-1}(\frac{c}{2})}  y\frac{1}{(1+y^2)^\frac{5}{4}} \frac{1}{\sqrt{s+y}} \; \mathrm{d}y. 
\end{align}
Now we use that $G^{-1}(\frac{c}{2}) = \frac{\theta s}{\theta+b}$ (cf. \eqref{eq:NewParameters})
to rewrite
\begin{equation}
    h =- \frac{s^\frac{1}{2}}{\sqrt{\theta+b}} \int_{-s}^{-\theta s/(\theta+b)}  y \frac{1}{(1+y^2)^\frac{5}{4}} \frac{1}{\sqrt{s+y}} \; \mathrm{d}y.
\end{equation}
The substitution $y= -s z$ yields 
\begin{equation}
    h = -\frac{s^\frac{1}{2}}{\sqrt{\theta+b}} \int_1^{\theta/(\theta+b)} (-s) \frac{-sz}{(1+s^2z^2)^\frac{5}{4}} \frac{1}{\sqrt{s-sz}} \; \mathrm{d}z = \frac{1}{\sqrt{\theta+b}} \int_{\theta/(\theta+b)}^1 \frac{s^2z}{(1+s^2z^2)^\frac{5}{4}} \frac{1}{\sqrt{1-z}} \; \mathrm{d}z. 
 \end{equation}
\end{proof}
\begin{corollary}\label{cor:sharpnonexistence} Suppose that $\theta > c_0^2$. 
    If $\mathcal{E}_\theta$ attains a minimizer in $ \mathcal{A}_{\rm sym}(h)$ then 
    \begin{equation}\label{eq:supremumh}
        h \leq \sup_{b,s > 0}  \frac{1}{\sqrt{\theta+b}}\int_{\theta/(\theta+b)}^1 \frac{s^2y}{(1+s^2y^2)^\frac{5}{4}} \frac{1}{\sqrt{1-y}} \; \mathrm{d}y =: h_0.
    \end{equation}
\end{corollary}
\begin{proof}
  Suppose that $\mathcal{E}$ attains a minimizer $u \in \mathcal{A}_{\rm sym}(h)$ . If $\{u = 0 \} = \emptyset$ we obtain $u> 0$ and thus 
   \begin{equation}
       \mathcal{E}(u) \geq \theta |\{ u > 0 \}|= \theta | (0,1)| = \theta > c_0^2,
   \end{equation}
   contradicting Proposition \ref{prop:theta}. If $\{ u = 0 \} \neq \emptyset$ we can apply Proposition \ref{Prop:4.1} and infer $h\leq  h_0$. 
\end{proof}
Next, we show that the quantity $h_0$ defined in the previous corollary is finite and also decreases with $\theta$. 
\begin{prop}\label{prop.secontpart1} Let $h_0$ be as in Corollary \ref{cor:sharpnonexistence}. Then
    \begin{equation}
        h_0 \leq \frac{1}{\sqrt{\theta}} \sup_{r > 0} \frac{2r}{\cosh(r)}.
    \end{equation}
\end{prop}
\begin{proof}
    For $s,b > 0$ we compute
    \begin{align}
         \frac{1}{\sqrt{\theta+b}}\int_{\theta/(\theta+b)}^1 \frac{s^2y}{(1+s^2y^2)^\frac{5}{4}} \frac{1}{\sqrt{1-y}} \; \mathrm{d}y & \leq \frac{1}{\sqrt{\theta+b}}\int_{\theta/(\theta+b)}^1 \frac{s^2y}{1+s^2y^2} \frac{1}{\sqrt{1-y}} \; \mathrm{d}y 
        \\ & \leq  \frac{1}{\sqrt{\theta+b}}\int_{\theta/(\theta+b)}^1 \frac{s^2y}{s^2y^2} \frac{1}{\sqrt{1-y}} \; \mathrm{d}y = \frac{1}{\sqrt{\theta+b}}\int_{\theta/(\theta+b)}^1  \frac{1}{y\sqrt{1-y}} \; \mathrm{d}y.
    \end{align}
    The substitution $u = \sqrt{1-y}$ yields 
    \begin{align}
         & \frac{1}{\sqrt{\theta+b}}\int_{\theta/(\theta+b)}^1 \frac{s^2y}{(1+s^2y^2)^\frac{5}{4}} \frac{1}{\sqrt{1-y}} \; \mathrm{d}y  \leq \frac{2}{\sqrt{\theta+b}}  \int_0^{\sqrt{\frac{b}{b+\theta}}} \frac{1}{1-u^2} \; \mathrm{d}u
         \\ & = \frac{1}{\sqrt{\theta+b}} \int_0^{\sqrt{\frac{b}{b+\theta}}} \frac{1}{1-u} + \frac{1}{1+u} \; \mathrm{d}u
          = \frac{1}{\sqrt{\theta+b}}  \log \frac{1+ \sqrt{\frac{b}{b+\theta}}}{1- \sqrt{\frac{b}{b+\theta}}} = \frac{1}{\sqrt{\theta+b}} \log \frac{\sqrt{b}+ \sqrt{b+\theta}}{\sqrt{b+\theta}-\sqrt{b}}
         \\ &= \frac{1}{\sqrt{\theta+b}} \log \left( \frac{(\sqrt{b}+ \sqrt{b+\theta})^2}{\theta}   \right) = \frac{2}{\sqrt{\theta+b}} \log\frac{\sqrt{b}+ \sqrt{b+\theta}}{\sqrt{\theta}}. 
    \end{align}
    Defining $b' := \theta b$ we obtain $\frac{2}{\sqrt{\theta}}  \frac{1}{\sqrt{1+b'}} \log( \sqrt{b'} + \sqrt{b' + 1})$. 
    Now define $x(b') := \sqrt{b'}$. Then clearly $\log(\sqrt{b'} + \sqrt{b'+1}) = \log (x(b')+ \sqrt{x(b')^2+1}) = \mathrm{Arsinh}(x(b'))$ and thus we have 
    \begin{equation}
         \frac{1}{\sqrt{\theta+b}}\int_{\theta/(\theta+b)}^1 \frac{s^2y}{(1+s^2y^2)^\frac{5}{4}} \frac{1}{\sqrt{1-y}} \; \mathrm{d}y \leq \frac{1}{\sqrt{\theta}} \frac{2\mathrm{Arsinh}(x(b'))}{\sqrt{1+x(b')^2}}.
    \end{equation}
    Defining $r = \mathrm{Arsinh}(x(b'))$ we obtain $\sqrt{1+x(b')^2} = \sqrt{1+ \sinh(r)^2} = \cosh(r)$ and thus 
    \begin{equation}
        \frac{1}{\sqrt{\theta+b}}\int_{\theta/(\theta+b)}^1 \frac{s^2y}{(1+s^2y^2)^\frac{5}{4}} \frac{1}{\sqrt{1-y}} \; \mathrm{d}y \leq \frac{1}{\sqrt{\theta}} \frac{2r}{\cosh(r)}.
    \end{equation}
    The claim follows. 
\end{proof}
\begin{remark}\label{rem:44}
    Note that $\sup\limits_{r > 0} \frac{2r}{\cosh(r)} = \max\limits_{r > 0} \frac{2r}{\cosh(r)}$ and the maximum is attained at some unique $r_*>0$. This is readily seen as follows: We have $\frac{d}{dr} \frac{2r}{\cosh(r)} = \frac{2(\cosh(r)- r \sinh(r))}{\cosh^2(r)}$. Now define $e: [0,\infty) \to \mathbb{R}, r \mapsto \cosh(r)- r \sinh(r)$. Given that $e(0)=1, \lim\limits_{r \to \infty} e(r) =- \infty $ and $e'(r) = - r \cosh(r) < 0$ for $r> 0$ we have that $e$ changes sign exactly once. This implies that there exists a unique $r_* > 0$ such that  $r \mapsto \frac{2r}{\cosh(r)}$ is strictly increasing  in $(0,r_*]$ and strictly decreasing in $[r_*,\infty)$. The claim follows.  
\end{remark}
\begin{lemma}\label{lem:secondpart2}
\begin{equation}
    \sup_{r> 0}  \frac{2r}{\cosh(r)} < \frac{8}{c_0^2}.
\end{equation}
\end{lemma}
\begin{proof}
    Define $f: (0, \infty) \to \mathbb{R}, r \mapsto 8\cosh(r) - 2c_0^2r$. Now, given that $f'(r) = 8\sinh(r) - 2c_0^2$ is increasing in $r$ and attains negative values in a neighborhood of $0$, we find that $f$ is decreasing until $r_0 := \mathrm{Arsinh}(\frac{c_0^2}{4})$ and increasing afterwards. Hence 
    \begin{equation}
        f(r) \geq f(r_0) = 8 \cosh(r_0) - 2c_0^2 r_0 = 8 \sqrt{1+ \sinh^2(r_0)} - 2 c_0^2 r_0 = 8  \sqrt{1+ (\tfrac{c_0^2}{4})^2} - 8 \tfrac{c_0^2}{4} \mathrm{Arsinh}(\tfrac{c_0^2}{4}).
    \end{equation}
    If we define $g(z) := 8( \sqrt{1+z^2} - z\mathrm{Arsinh}(z))$ then we have $f(r_0)= g(\frac{c_0^2}{4})$.
    Clearly $g'(z) = -8\mathrm{Arsinh}(z) <0$, which means that $g$ is strictly decreasing. Using this and that $c_0^2 < 6$ (which we prove independently in Lemma \ref{lem:cnullsquaresmaller6})  we obtain 
    $f(r_0) =  g(\tfrac{c_0^2}{4}) > g(\tfrac{3}{2}) = 8 (\tfrac{\sqrt{13}}{2} - \tfrac{3}{2} \log(\tfrac{3}{2} + \tfrac{\sqrt{13}}{2})) = 4 ( \sqrt{13}- 3 \log(\tfrac{3+ \sqrt{13}}{2})) =: \alpha.$
    We claim that the expression  $\alpha$ is positive. Indeed, this is equivalent to $\log(\frac{3+\sqrt{13}}{2}) < \frac{\sqrt{13}}{3}$ or alternatively $e^\frac{\sqrt{13}}{3} > \frac{3+\sqrt{13}}{2}$. This can be readily verified by Taylor expansion. Indeed,
    \begin{equation}
        e^\frac{\sqrt{13}}{3}- \frac{3}{2}- \frac{\sqrt{13}}{2} > \sum_{k = 0}^3 \frac{1}{k!} \left( \frac{\sqrt{13}}{3} \right)^k - \frac{3}{2}- \frac{\sqrt{13}}{2} = 1+ \frac{\sqrt{13}}{3}+ \frac{13}{18} + \frac{13}{162} \sqrt{13}  - \frac{3}{2}- \frac{\sqrt{13}}{2} = \frac{4\sqrt{13}-9}{81},
    \end{equation}
    which is positive as $(4 \sqrt{13})^2 = 16\cdot 13 > 100 > 81=  9^2.$ We have now shown that $f(r) \geq \alpha$ for all $r \in (0,\infty)$. This can be rearranged as 
    $
        \frac{2r}{\cosh(r)} \leq \frac{8}{c_0^2}- \frac{\alpha}{\cosh(r)},
    $
    which now holds for all $r> 0$. As a result we have, with $r_*$ as in  Remark \ref{rem:44},
   $
        \frac{2r}{\cosh(r)} \leq \frac{2r_*}{\cosh(r_*)} \leq  \frac{8}{c_0^2}- \frac{\alpha}{\cosh(r_*)}.
 $
    Hence, $\sup\limits_{r> 0} \frac{2r}{\cosh(r)} \leq  \frac{8}{c_0^2}- \frac{\alpha}{\cosh(r_*)}< \frac{8}{c_0^2}$. The claim follows. 
\end{proof}
To complete the previous proof we derive the following lemma.
\begin{lemma}\label{lem:cnullsquaresmaller6}
    $c_0^2 < 6.$ 
\end{lemma}
\begin{proof}
    The substitution $s= \tan(t)$ implies $$c_0 = \int_{\mathbb{R}} \frac{1}{(1+s^2)^\frac{5}{4}} \; \mathrm{d}s  = \int_{-\frac{\pi}{2}}^\frac{\pi}{2} \frac{1}{(1+\tan^2(t))^\frac{1}{4}} \; \mathrm{d}t = \int_{-\pi/2}^{\pi/2} \sqrt{\cos(t)} \; \mathrm{d}t = 2 \int_{0}^{\pi/2} \sqrt{\cos(t)} \; \mathrm{d}t.$$ 
Using $\cos(t) = 1- 2 \sin^2(\frac{t}{2})$ and expanding $\sqrt{1- 2 \sin^2(\frac{t}{2})}$  with the Taylor series of $x \mapsto \sqrt{1+x}$ yields 
\begin{align}\label{eq:cnought}
    \int_0^\frac{\pi}{2} \sqrt{\cos(t)} \; \mathrm{d}t = \int_0^\frac{\pi}{2} \sum_{k = 0}^\infty \binom{\frac{1}{2}}{k} 2^k (-1)^k \sin^{2k}(\tfrac{t}{2}) \; \mathrm{d}t   = \int_0^\frac{\pi}{2}  \left( 1 - \sum_{k = 1}^\infty \binom{\frac{1}{2}}{k} 2^k (-1)^{k-1} \sin^{2k}(\tfrac{t}{2}) \right) \; \mathrm{d}t. 
\end{align}
Note that for all $k \geq 1$ one has $\binom{\frac{1}{2}}{k} (-1)^{k-1} \geq 0$ (which is easily shown by induction). 
Hence we may estimate 
\begin{equation}
    \int_0^\frac{\pi}{2} \sqrt{\cos(t)} \; \mathrm{d}t \leq \int_0^\frac{\pi}{2}  \left( 1 - \sum_{k = 1}^3 \binom{\frac{1}{2}}{k} 2^k (-1)^{k-1} \sin^{2k}(\tfrac{t}{2}) \right) \; \mathrm{d}t = \frac{\pi}{2} - \sum_{k = 1}^3 \binom{\frac{1}{2}}{k} 2^k (-1)^{k-1} \int_0^\frac{\pi}{2} \sin^{2k}(\tfrac{t}{2}) \; \mathrm{d}t. 
\end{equation}
Now we use that 
$
    \int_0^\frac{\pi}{2} \sin^2(\tfrac{t}{2}) \; \mathrm{dt} = \tfrac{\pi}{4} - \tfrac{1}{2}, \;  \int_0^\frac{\pi}{2} \sin^4(\tfrac{t}{2}) \; \mathrm{dt} = \tfrac{3\pi}{16} - \tfrac{1}{2},  \; \int_0^\frac{\pi}{2} \sin^6(\tfrac{t}{2}) \; \mathrm{dt}= \tfrac{5\pi}{32}- \tfrac{11}{24}
$
to obtain 
\begin{equation}
    \int_0^\frac{\pi}{2} \sqrt{\cos(t)} \; \mathrm{d}t \leq \tfrac{\pi}{2} - \left(\tfrac{\pi}{4}-\tfrac{1}{2} \right) - \tfrac{1}{2} \left(  \tfrac{3\pi}{16} - \tfrac{1}{2} \right) - \tfrac{1}{2} \left( \tfrac{5\pi}{32}- \tfrac{11}{24} \right)  = \tfrac{5}{64}\pi + \tfrac{47}{48}.  
\end{equation}
This implies with \eqref{eq:cnought} that $c_0 \leq \frac{5}{32}\pi+ \frac{47}{24}.$
One can use next that $\pi < \frac{22}{7}$ (cf. \cite{Dalzell}) and estimate this quantity by $\frac{823}{336}$. Hence $c_0^2  \leq \frac{823^2}{336^2} =\frac{677329}{112896} < 6$.

\end{proof}
\begin{proof}[Proof of Theorem \ref{thm2}]
    The first sentence of the claim follows directly from Corollary \ref{cor:sharpnonexistence}. The second sentence follows from Proposition \ref{prop.secontpart1} and Lemma \ref{lem:secondpart2}.
\end{proof}

\begin{remark}
    We have to leave the question open whether the value in Theorem \ref{thm2} is an optimal choice for $h_*$. This shall be subject of future research.  
\end{remark}

\bibliographystyle{abbrv}
\bibliography{ref_Mueller-Yoshizawa}

\begin{bibdiv}
\begin{biblist}

\bib{Bartels}{article}{
      author={Bartels, S\"oren},
      author={Weyer, Pascal},
       title={Computing confined elasticae},
        date={2022},
        ISSN={2731-4235},
     journal={Adv. Contin. Discrete Models},
       pages={Paper No. 58, 16},
         url={https://doi.org/10.1186/s13662-022-03731-7},
      review={\MR{4503111}},
}

\bib{DDObstacle}{article}{
      author={Dall'Acqua, Anna},
      author={Deckelnick, Klaus},
       title={An obstacle problem for elastic graphs},
        date={2018},
        ISSN={0036-1410,1095-7154},
     journal={SIAM J. Math. Anal.},
      volume={50},
      number={1},
       pages={119\ndash 137},
         url={https://doi.org/10.1137/17M111701X},
      review={\MR{3742685}},
}

\bib{ObstacleP}{article}{
      author={Dall'Acqua, Anna},
      author={M\"uller, Marius},
      author={Okabe, Shinya},
      author={Yoshizawa, Kensuke},
       title={An obstacle problem for the {$p$}-elastic energy},
        date={2024},
        ISSN={0944-2669,1432-0835},
     journal={Calc. Var. Partial Differential Equations},
      volume={63},
      number={6},
       pages={Paper No. 145, 43},
         url={https://doi.org/10.1007/s00526-024-02752-2},
      review={\MR{4761861}},
}

\bib{Network1}{article}{
      author={Dall'Acqua, Anna},
      author={Novaga, Matteo},
      author={Pluda, Alessandra},
       title={Minimal elastic networks},
        date={2020},
        ISSN={0022-2518,1943-5258},
     journal={Indiana Univ. Math. J.},
      volume={69},
      number={6},
       pages={1909\ndash 1932},
         url={https://doi.org/10.1512/iumj.2020.69.8036},
      review={\MR{4170083}},
}

\bib{Network2}{article}{
      author={Dall'Acqua, Anna},
      author={Pluda, Alessandra},
       title={Some minimization problems for planar networks of elastic curves},
        date={2017},
        ISSN={2353-3382},
     journal={Geom. Flows},
      volume={2},
      number={1},
       pages={105\ndash 124},
         url={https://doi.org/10.1515/geofl-2017-0005},
      review={\MR{3745454}},
}

\bib{Dalzell}{article}{
      author={Dalzell, D.~P.},
       title={On {$22/7$}},
        date={1944},
        ISSN={0024-6107,1469-7750},
     journal={J. London Math. Soc.},
      volume={19},
       pages={133\ndash 134},
         url={https://doi.org/10.1112/jlms/19.75_Part_3.133},
      review={\MR{13425}},
}

\bib{DMN18}{article}{
      author={Dayrens, Fran\c~cois},
      author={Masnou, Simon},
      author={Novaga, Matteo},
       title={Existence, regularity and structure of confined elasticae},
        date={2018},
        ISSN={1292-8119,1262-3377},
     journal={ESAIM Control Optim. Calc. Var.},
      volume={24},
      number={1},
       pages={25\ndash 43},
         url={https://doi.org/10.1051/cocv/2016073},
      review={\MR{3764132}},
}

\bib{Network3}{article}{
      author={Del~Nin, Giacomo},
      author={Pluda, Alessandra},
      author={Pozzetta, Marco},
       title={Degenerate elastic networks},
        date={2021},
        ISSN={1050-6926,1559-002X},
     journal={J. Geom. Anal.},
      volume={31},
      number={6},
       pages={6128\ndash 6170},
         url={https://doi.org/10.1007/s12220-020-00521-z},
      review={\MR{4267640}},
}

\bib{DKV1}{article}{
      author={Dipierro, Serena},
      author={Karakhanyan, Aram},
      author={Valdinoci, Enrico},
       title={A free boundary problem driven by the biharmonic operator},
        date={2020},
        ISSN={2578-5885,2578-5893},
     journal={Pure Appl. Anal.},
      volume={2},
      number={4},
       pages={875\ndash 942},
         url={https://doi.org/10.2140/paa.2020.2.875},
      review={\MR{4229229}},
}

\bib{DKV2}{article}{
      author={Dipierro, Serena},
      author={Karakhanyan, Aram~L.},
      author={Valdinoci, Enrico},
       title={Limit behaviour of a singular perturbation problem for the biharmonic operator},
        date={2019},
        ISSN={0095-4616,1432-0606},
     journal={Appl. Math. Optim.},
      volume={80},
      number={3},
       pages={679\ndash 713},
         url={https://doi.org/10.1007/s00245-019-09598-7},
      review={\MR{4026596}},
}

\bib{EVGARIEPY2015}{book}{
      author={Evans, Lawrence~C.},
      author={Gariepy, Ronald~F.},
       title={Measure theory and fine properties of functions},
     edition={Revised},
      series={Textbooks in Mathematics},
   publisher={CRC Press, Boca Raton, FL},
        date={2015},
        ISBN={978-1-4822-4238-6},
      review={\MR{3409135}},
}

\bib{HG1}{article}{
      author={Grunau, Hans-Christoph},
      author={M\"uller, Marius},
       title={A biharmonic analogue of the {A}lt-{C}affarelli problem},
        date={2024},
        ISSN={0025-5831,1432-1807},
     journal={Math. Ann.},
      volume={390},
      number={4},
       pages={5259\ndash 5297},
         url={https://doi.org/10.1007/s00208-024-02883-z},
      review={\MR{4816111}},
}

\bib{HG2}{article}{
      author={Grunau, Hans-Christoph},
      author={M\"uller, Marius},
       title={Global minimizers for a two-sided biharmonic alt-caffarelli problem},
        date={2026},
     journal={SIAM J. Math. Anal. (accepted)},
}

\bib{GrunauOkabe1}{article}{
      author={Grunau, Hans-Christoph},
      author={Okabe, Shinya},
       title={Willmore obstacle problems under {D}irichlet boundary conditions},
        date={2023},
        ISSN={0391-173X,2036-2145},
     journal={Ann. Sc. Norm. Super. Pisa Cl. Sci. (5)},
      volume={24},
      number={3},
       pages={1415\ndash 1462},
         url={https://doi.org/10.2422/2036-2145.202105_064},
      review={\MR{4675964}},
}

\bib{GrunauOkabe2}{article}{
      author={Grunau, Hans-Christoph},
      author={Okabe, Shinya},
       title={Optimality of smallness conditions in {W}illmore obstacle problems under {D}irichlet boundary conditions},
        date={2025},
        ISSN={1468-1218,1878-5719},
     journal={Nonlinear Anal. Real World Appl.},
      volume={85},
       pages={Paper No. 104363, 10},
         url={https://doi.org/10.1016/j.nonrwa.2025.104363},
      review={\MR{4878867}},
}

\bib{NahonLamboley}{article}{
      author={Lamboley, Jimmy},
      author={Mickael, Nahon},
       title={Boundary regularity of a fourth order {A}lt-{C}affarelli problem and applications to the minimization of the critical buckling load},
        date={2025},
     journal={Preprint},
}

\bib{Miura16}{article}{
      author={Miura, Tatsuya},
       title={Singular perturbation by bending for an adhesive obstacle problem},
        date={2016},
        ISSN={0944-2669,1432-0835},
     journal={Calc. Var. Partial Differential Equations},
      volume={55},
      number={1},
       pages={Art. 19, 24},
         url={https://doi.org/10.1007/s00526-015-0941-z},
      review={\MR{3456944}},
}

\bib{Miura21}{article}{
      author={Miura, Tatsuya},
       title={Polar tangential angles and free elasticae},
        date={2021},
     journal={Math. Eng.},
      volume={3},
      number={4},
       pages={Paper No. 034, 12},
         url={https://doi.org/10.3934/mine.2021034},
      review={\MR{4147578}},
}

\bib{ObstacleProbl}{article}{
      author={M\"uller, Marius},
       title={An obstacle problem for elastic curves: existence results},
        date={2019},
        ISSN={1463-9963,1463-9971},
     journal={Interfaces Free Bound.},
      volume={21},
      number={1},
       pages={87\ndash 129},
         url={https://doi.org/10.4171/IFB/418},
      review={\MR{3951579}},
}

\bib{MariusAMOP}{article}{
      author={M\"uller, Marius},
       title={The elastic flow with obstacles: small obstacle results},
        date={2021},
        ISSN={0095-4616,1432-0606},
     journal={Appl. Math. Optim.},
      volume={84},
       pages={S355\ndash S402},
         url={https://doi.org/10.1007/s00245-021-09773-9},
      review={\MR{4316789}},
}

\bib{MariusAMPA}{article}{
      author={M\"uller, Marius},
       title={The biharmonic {A}lt-{C}affarelli problem in 2{D}},
        date={2022},
        ISSN={0373-3114,1618-1891},
     journal={Ann. Mat. Pura Appl. (4)},
      volume={201},
      number={4},
       pages={1753\ndash 1799},
         url={https://doi.org/10.1007/s10231-021-01178-3},
      review={\MR{4454382}},
}

\bib{MariusIFB}{article}{
      author={M\"uller, Marius},
       title={Polyharmonic equations involving surface measures},
        date={2024},
        ISSN={1463-9963,1463-9971},
     journal={Interfaces Free Bound.},
      volume={26},
      number={1},
       pages={61\ndash 78},
         url={https://doi.org/10.4171/ifb/503},
      review={\MR{4705642}},
}

\bib{MariusCVPD}{article}{
      author={M\"uller, Marius},
       title={An anisotropic {A}lt-{C}affarelli problem of higher order},
        date={2025},
     journal={Preprint},
}

\bib{RoegerRupp}{article}{
      author={Röger, Matthias},
      author={Rupp, Fabian},
       title={Elastic bilayer membranes under confinement -- existence, regularity, and rigidity},
        date={2025},
     journal={Preprint},
}

\bib{KY2021}{article}{
      author={Yoshizawa, Kensuke},
       title={A remark on elastic graphs with the symmetric cone obstacle},
        date={2021},
        ISSN={0036-1410,1095-7154},
     journal={SIAM J. Math. Anal.},
      volume={53},
      number={2},
       pages={1857\ndash 1885},
         url={https://doi.org/10.1137/19M1307901},
      review={\MR{4238008}},
}

\end{biblist}
\end{bibdiv}

\end{document}